\documentclass{amsart}

\usepackage{amscd,amssymb}
\usepackage[all]{xy}

\usepackage{hyperref}

\def\Z{{\mathbb Z}}
\def\Q{{\mathbb Q}}
\def\R{{\mathbb R}}
\def\C{{\mathbb C}}
\def\P{{\mathbb P}}

\def\bH{{\mathbb H}}

\def\A{{\mathcal A}}
\def\H{{\mathcal H}}
\def\K{{\mathcal K}}
\def\M{{\mathcal M}}
\def\O{{\mathcal O}}

\def\T{{\mathcal T}}
\def\U{{\mathcal U}}
\def\X{{\mathcal X}}
\def\cC{{\mathcal C}}

\def\cG{{\mathcal G}}
\def\cI{{\mathcal I}}
\def\cP{{\mathcal P}}

\def\a{{\mathfrak a}}

\def\k{{\mathfrak k}}
\def\n{{\mathfrak n}}
\def\p{{\mathfrak p}}

\def\t{{\mathfrak t}}
\def\u{{\mathfrak u}}

\def\G{{\Gamma}}
\def\La{\Lambda}

\def\Ghat{{\hat{\G}}}

\def\Qbar{{\overline{\Q}}}

\def\Xbar{\overline{X}}

\def\rhotilde{{\tilde{\rho}}}
\def\rhohat{{\hat{\rho}}}

\def\cXtilde{{\widetilde{\X}}}

\def\l{{\ell}}
\def\v{{\overrightarrow{01}}}
\def\w{{\omega}}

\def\arith{\mathrm{arith}}

\def\prol{{(\ell)}}

\def\Vec{\mathrm{Vec}}
\def\real{\mathrm{real}}

\def\Ql{{\Q_\ell}}
\def\Qlx{{\Q_\ell^\times}}
\def\Zl{{\Z_\ell}}
\def\Zlx{{\Z_\ell^\times}}

\def\Gl{{G_\ell}}

\def\XFS{{X_{F,S}}}
\def\GFS{{G_{F,S}}}
\def\OFS{{\O_{F,S}}}

\def\AFSl{{\A_{F,S}^\l}}
\def\KFSl{{\K_{F,S}^\l}}
\def\aFSl{{\a_{F,S}^\l}}
\def\kFSl{{\k_{F,S}^\l}}

\def\MTMlFS{{\T_\l(\XFS)}}

\def\I{{I_\l}}
\def\cIell{{\mathcal{I}_\l}}

\def\Pminus{{\P^1 - \{0,1,\infty\}}}
\def\PminusC{{\P^1(\C) - \{0,1,\infty\}}}

\def\Hcts{H_{\mathrm{cts}}}
\def\Hctsf{H_{\mathrm{cts}f}}
\def\Het{H_\et}
\def\Hmot{H_{\mathrm{mot}}}
\def\Hdel{H_{\mathcal D}}
\def\Hc{H_{\mathcal C}}
\def\cts{\mathrm{cts}}
\def\et{\mathrm{\acute{e}t}}

\def\Gm{{\mathbb{G}_m}}

\def\GSp{{\mathrm{GSp}}}
\def\Sp{{\mathrm{Sp}}}
\def\GL{{\mathrm{GL}}}

\def\dot{{\bullet}}

\def\blank{\phantom{x}}

\def\mot{{\mathrm{mot}}}
\def\limproj#1{\lim_{\stackrel{\longleftarrow}{#1}}}
\def\liminj#1{\lim_{\stackrel{\longrightarrow}{#1}}}

\newcommand\Li{\operatorname{Li}}
\newcommand\im{\operatorname{im}}
\newcommand\ord{\operatorname{ord}}

\newcommand\Spec{\operatorname{Spec}}
\newcommand\Hom{\operatorname{Hom}}

\newcommand\Ext{\operatorname{Ext}}
\newcommand\sExt{\operatorname{\mathcal E\!\mathit{xt}}}
\newcommand\Aut{\operatorname{Aut}}

\newcommand\Out{\operatorname{Out}}

\newcommand\OutDer{\operatorname{OutDer}}
\newcommand\Gal{\operatorname{Gal}}
\newcommand\Gr{\operatorname{Gr}}

\newtheorem{theorem}{Theorem}[section]
\newtheorem{lemma}[theorem]{Lemma}
\newtheorem{proposition}[theorem]{Proposition}
\newtheorem{corollary}[theorem]{Corollary}

\theoremstyle{definition}
\newtheorem{definition}[theorem]{Definition}
\newtheorem{example}[theorem]{Example}
\newtheorem{conjecture}[theorem]{Conjecture}
\newtheorem{postulate}[theorem]{Postulate}

\theoremstyle{remark}
\newtheorem{remark}[theorem]{Remark}
\begin{document}

\title{Tannakian Fundamental Groups Associated to Galois Groups}

\author{Richard Hain}
\address{Department of Mathematics\\ Duke University\\
Durham, NC 27708-0320}
\email{hain@math.duke.edu}
\thanks{The first author was supported in part by grants from the National
Science Foundation. The second author was supported in part by a Mombusho
Grant and also by MSRI during a visit in the fall of 1999.}

\author{Makoto Matsumoto}
\address{Department of Mathematics\\ 
Faculty of Science \\ Hiroshima University\\
Higashi-Hiroshima, 739-8526 Japan}
\email{m-mat@math.sci.hiroshima-u.ac.jp}

%\date{\today}

%\subjclass{Primary xxxxx; Secondary xxxxx}

%\keywords{}

\maketitle

\section{Introduction}

The goal of this paper is to provide background, heuristics and motivation for
several conjectures of Deligne \cite[8.2, p.~163]{deligne}, \cite[8.9.5,
p.~168]{deligne} and \cite[p.~300]{Ihara3} and Goncharov
\cite[Conj.~2.1]{goncharov}, presumably along the lines used to arrive at them.
A complete proof of the third of these conjectures, and partial solutions of
the remaining three are given in \cite{hain-matsumoto}.\footnote{After writing
that paper, we learned from Goncharov that proofs of $\l$-adic versions of
\cite[8.2, p.~163]{deligne} and \cite[8.9.5, p.~168]{deligne} had previously
been given in the unpublished manuscript \cite{beilinson-deligne:orig} of
Beilinson and Deligne.} A second goal of this paper is to show that the
weighted completion of a profinite group, developed in \cite{hain-matsumoto},
and a key ingredient in the proofs referred to above, can be defined as the
tannakian fundamental group\footnote{A tannakian category $\T$ with fiber
functor $\w$ is equivalent to the category of representations of the
automorphism group of $\w$. We shall refer to this proalgebraic group as the
{\it tannakian fundamental group} of $\T$ with respect to the base point $\w$.
Basic material on tannakian categories, such as their definition, can be found
in \cite{dmos}.} of certain categories of modules of the group. This should
help clarify the role of weighted completion in \cite{hain-matsumoto}.

\section{Motivic Cohomology}

It is believed that there is a universal cohomology theory, called
{\it motivic cohomology}. It should be defined for all schemes $X$. It is
indexed by two integers $m$ and $n$. The coefficient
ring $\La$ is typically $\Z$, $\Z/N$, $\Zl$, $\Q$ or $\Ql$; the
corresponding motivic cohomology group is denoted
$$
\Hmot^m(X,\La(n)).
$$
There should be cup products
\begin{equation}
\label{cup_prod}
\Hmot^{m_1}(X,\La(n_1)) \otimes \Hmot^{m_2}(X,\La(n_2))
\to \Hmot^{m_1+m_2}(X,\La(n_1+n_2)).
\end{equation}
Motivic cohomology should have the following universal mapping property:
if $\Hc^\dot(\blank,\Lambda(\blank))$ is any Bloch-Ogus cohomology theory
\cite{bloch-ogus} (such as \'etale cohomology, Deligne cohomology, Betti
(i.e., singular) cohomology, de Rham cohomology, and crystalline cohomology)
there should be a unique natural transformation
$$
\Hmot^m(\blank,\La(n)) \to \Hc^m(\blank,\La(n))
$$
compatible with products and Chern classes
$$
c_n : K_m(X) \to \Hc^{2n-m}(X,\La(n)),
$$
where $K_\dot$ denotes Quillen's algebraic $K$-group functor, \cite{quillen}.

\subsection*{Beilinson's definition}

Beilinson \cite{beilinson} observed that the motivic cohomology with $\Q$
coefficients of a large class of schemes could be defined in terms of
Quillen's algebraic $K$-theory \cite{quillen}.

Suppose that $X$ is the spectrum of the ring of $S$-integers in a number
field or a smooth scheme over a perfect field. Denote the algebraic $K$-theory
of $X$ by $K_\dot(X)$. As in the case of topological
$K$-theory, there are Adams operations (\cite{hiller}, \cite{kratzer})
$$
\psi^k : K_\dot(X) \to K_\dot(X),
$$
defined for all $k \in \Z_+$. They can be simultaneously diagonalized:
$$
K_m(X)\otimes \Q = \bigoplus_{n\in \Z} K_m(X)^{(n)}
$$
where $\psi^k$ acts as $k^n$ on $K_m(X)^{(n)}$.

\begin{definition}[Beilinson]
\label{def}
For a ring $\La$ containing $\Q$, define the motivic cohomology groups
of $X$ by
$$
\Hmot^m(X,\La(n)) = K_{2n-m}(X)^{(n)}\otimes_\Q \La.
$$
\end{definition} 

The ring structure of $K_\dot(X)$ induces a cup product (\ref{cup_prod}) as
$$
\psi^k(xy) = \psi^k(x)\psi^k(y)\quad x,y \in K_\dot(X).
$$
Motivation for Beilinson's definition comes from topological $K$-theory and
can be found in the introduction of \cite{bms}.

If $X$ is smooth, then it follows from a result of Grothendieck (see
\cite{borel-serre}) that
$$
\Hmot^{2n}(X,\Q(n)) \cong CH^n(X)\otimes\Q.
$$
In the next section, we present the well-known computation of the motivic
cohomology of the ring of $S$-integers in a number field.

\begin{proposition}
There are Chern classes
$$
c^\mot_j : K_m(X) \to \Hmot^{2j-m}(X,\Q(j))
$$
such that for each Bloch-Ogus cohomology theory $\Hc^\dot(\blank,\La(\blank))$,
where $\La$ contains $\Q$, there is a natural transformation
$$
\Hmot^\dot(\blank,\La(\blank)) \to \Hc^\dot(\blank,\La(\blank))
$$
that is compatible with Chern classes.
\end{proposition}

\begin{proof}
The basic tool needed to construct the natural transformations to other 
cohomology theories is the theory of Chern classes
$$
c_j : K_m(X) \to \Hc^{2j-m}(X,\Z(j))
$$
constructed by Beilinson \cite{beilinson} and Gillet \cite{gillet} for a
very large set of cohomology theories $\Hc^\dot$ that includes all Bloch-Ogus
cohomology theories. These give rise to the Chern character maps
$$
ch : K_m(X) \to \prod_{j\ge 0} \Hc^{2j-m}(X,\Q(j)).
$$
The degree $j$ part of this
$$
ch_j : K_m(X) \to \Hc^{2j-m}(X,\Q(j))
$$
is a homogeneous polynomial of degree $j$ in the Chern classes, just
as in the topological case. The key point is the compatibility with
the Adams operations which implies that the restriction of $ch_k$
to $K_m(X)^{(j)}$ vanishes unless $k=j$. It follows that $ch_j$ factors
through the projection onto $K_m(X)^{(j)}$:
$$
\xymatrix{
K_m(X) \ar[r]^{ch_j\qquad} \ar[d]_{\mathrm{proj}} & \Hc^{2j-m}(X,\Q(j))\cr
K_m(X)^{(j)} \ar[ur]
}
$$
Thus the Chern character induces a natural transformation
$$
\Hmot^m(X,\Q(n)) \to \Hc^m(X,\Q(n)).
$$
It is a ring homomorphism as the Chern character is. Define
$$
ch^\mot_j : K_m(X) \to \Hmot^{2j-m}(X,\Q(j))
$$
to be the projection $K_m(X) \to K_m(X)^{(j)}$. From this, one can inductively
construct Chern classes $c^\mot_j : K_m(X) \to \Hmot^{2j-m}(X,\Q(j))$.
Compatibility with the Chern classes $c_j : K_m(X) \to \Hc^{2j-m}(X,\La(j))$ is
automatic and guarantees the uniqueness of natural transformation $\Hmot^\dot
\to \Hc^\dot$.
\end{proof}

\subsection*{The quest for cochains}

Beilinson's definition raises many questions and problems such as:
\begin{enumerate}
\item How does one define motivic cohomology with {\em integral} coefficients?
\item Find natural cochain complexes (a.k.a., motivic complexes) whose
homology groups are motivic cohomology.
\item Compute motivic cohomology groups.
\end{enumerate}

Bloch's higher Chow groups \cite{bloch} provide an integral version of
motivic cohomology as well as a chain complex whose homology is motivic
cohomology. (See also \cite{bloch:corr} and \cite{levine}.) One difficulty
with this approach is that, being based on algebraic cycles and rational
equivalence, it is difficult to compute with.

More fundamentally, one would also like motivic cohomology groups to be the ext
or hyper-ext groups associated to a suitable category of motives. In the ideal
case, this category will be tannakian after tensoring all its objects with
$\Q$, so that the category of $\Q$-motives will be equivalent to the category
of representations of a proalgebraic group defined over $\Q$. These goals have
been achieved to some degree. For all fields $k$, Voevodsky \cite{voevodsky}
and Levine \cite{levine:book} have each constructed a triangulated tensor
category of ``mixed motives over $k$''. For each scheme $X$, smooth and
quasi-projective over $k$, there is an object $M(X)$ in this category such that
$\Ext^\dot(\Z(-n),M(X))$ is isomorphic to the integral motivic cohomology
groups of $X$ (i.e., Bloch's higher Chow groups). However, the categories
obtained from the categories of Levine and Voevodsky by tensoring their objects
with $\Q$ are not tannakian.

One can also propose that there should be a tannakian category of {\it mixed
Tate motives} over a field $k$. The motivic cohomology of $k$ should be an ext
in this category. In the case where $k$ is a number field (or any field
satisfying Beilinson-Soul\'e vanishing), Levine \cite{levine:tate} has
constructed such a tannakian category of mixed Tate motives. Goncharov
\cite[p.~611]{goncharov:mtm} later proved a result similar to Levine's and
proved, in addition, that the bounded derived category of this tannakian
category of mixed Tate motives is equivalent to the full subcategory of mixed
Tate motives of the category of mixed motives over $k$.

An older and less fundamental approach to constructing categories of motives,
proposed by Deligne \cite{deligne} and Jannsen \cite{jannsen}, is to view them
as ``compatible systems of realizations''. These also form a tannakian
category. We shall take this approach in this paper as it is more accessible
and is more consistent with our point of view.

\section{The Motivic Cohomology of the Spectrum of a Ring of $S$-integers}

Basic results of Quillen \cite{quillen:fg} and Borel \cite{borel} give the
computation of the motivic cohomology of the spectra of rings of $S$-integers
in number fields. Suppose that $F$ is a number field with ring of integers
$\O_F$ and that $S$ is a finite subset of $\Spec \O_F$. Set $\XFS =
\Spec \O_F - S$. Set
\begin{equation}
\label{dn}
d_n = \ord_{s=1-n}\zeta_F(s) =
\begin{cases}
r_1 + r_2 - 1 & \text{ when } n=1,\cr
r_1 + r_2 & \text{ when $n$ is odd and } n > 1, \cr
r_2 & \text{ when $n$ is even,}
\end{cases}
\end{equation}
where $\zeta_F(s)$ denotes the Dedekind zeta function of $F$ and $r_1$,
and $r_2$ denote the number of real and complex places of $F$, respectively.

\begin{theorem}
\label{mot}
For all $n$ and $m$, $\Hmot^m(\XFS,\Q(n))$ is a finite dimensional rational
vector space whose dimension is given by
$$
\dim \Hmot^m(\XFS,\Q(n)) =
\begin{cases}
d_1 + \#S & m = n = 1, \cr
d_n & m=1 \text{ and } n > 1, \cr
0 & \text{otherwise.}
\end{cases}
$$
\end{theorem}

\begin{proof}
First suppose that $S$ is empty. Quillen \cite{quillen:fg} showed that each
$K$-group $K_m(\XFS)$ is a finitely generated abelian group. It follows that
each of the groups $\Hmot^j(\XFS,\Q(n))$ is finite dimensional. The rank of
$K_0(\XFS)$ is 1 and the rank of $K_1(\XFS)$ is $r_1 + r_2 - 1$ by the
Dirichlet Unit Theorem. The ranks of the remaining $K_m(\XFS)$ were computed
by Borel \cite{borel}. It is zero when $m$ is even and $>0$, and $d_n$
when $m = 2n-1 > 1$. It is easy to see that
$$
\Hmot^0(\XFS,\Q(0)) = K_0(\XFS) \otimes \Q \cong \Q.
$$
Borel \cite{borel:reg} constructed regulator mappings
$$
K_{2n-1}(\XFS) \to \R^{d_n}, \quad n > 0,
$$
and showed that each is injective mod torsion. Beilinson \cite{beilinson}
showed that Borel's regulator is a non-zero rational multiple of the regulator
mapping
$$
ch_n : K_{2n-1}(\XFS) \to \Hdel^1(\XFS,\R(n)) \cong \R^{d_n}
$$
to Deligne cohomology. The properties of the Chern character and Borel's
injectivity together imply that
$$
\Hmot^1(\XFS,\Q(n)) = K_{2n-1}(\XFS)^{(n)} = K_{2n-1}(\XFS)\otimes \Q
$$
and that $\Hmot^m(\XFS,\Q(n))$ vanishes when $m > 1$, and when $m=0$ and $n\neq
0$. The result when $S$ is non-empty follows from this using the localization
sequence \cite{quillen}, and the fact, due to Quillen \cite{quillen:fte}, that
the $K$-groups of finite fields are torsion groups in positive degree. Together
these imply that each prime removed adds one to the rank of $K_1$ and does not
change the rank of any other $K$-group.
\end{proof}

Denote the Galois group of the maximal algebraic extension of $F$, unramified
outside $S$, by $\GFS$. In this paper, a finite dimensional $\GFS$-module means
a finite dimensional $\Ql$-vector space with continuous $\GFS$-action. Denote
the category of $\Q$ mixed Hodge structures by $\H$.  Denote the ext functor in
the category of finite dimensional $\GFS$-modules by $\Ext_\GFS$, and the ext
functor in $\H$ by $\Ext_\H$. The results on regulators of Borel
\cite{borel:reg} and Soul\'e \cite{soule} can be stated as follows.

\begin{theorem}
\label{regulators}
The natural transformation from motivic to \'etale cohomology induces
isomorphisms
$$
\Hmot^1(\XFS,\Ql(n)) \cong \Het^1(\XFS,\Ql(n)) \cong \Ext^1_\GFS(\Ql,\Ql(n))
$$
for all $n\geq 1$.
The natural transformation from motivic to Deligne cohomology induces
injections
$$
\Hmot^1(\XFS,\Q(n)) \hookrightarrow \Hdel^1(\XFS,\Q(n)) =
\bigg[\bigoplus_{\nu : F \hookrightarrow \C}
\Ext^1_\H(\Q,\Q(n))\bigg]^{\Gal(\C/\R)}. \qed
$$
\end{theorem}

Thus each element $x$ of $K_{2n-1}(\XFS)$ determines an extension
$$
0 \to \Ql(n) \to E_{\l,x} \to \Ql(0) \to 0
$$
of $\l$-adic local systems over $\XFS$ and a $\Gal(\C/\R)$-equivariant
extension
$$
0 \to \Q(n) \to E_{\mathrm{Hodge},x} \to \Q(0) \to 0
$$
of mixed Hodge structures over $\XFS \otimes \C$. One can think of these
as the \'etale and Hodge realizations of $x\in K_{2n-1}(\XFS)$.

\section{Mixed Tate Motives}\label{sec:mtm}

As mentioned earlier, one approach to motivic cohomology is to postulate that
to each sufficiently nice scheme $X$ (say,  smooth and quasi-projective over a
field, or regular over a ring $\OFS$ of $S$-integers in a number field) one can
associate a category $\T(X)$ of {\it mixed Tate motives over} $X$. This should
satisfy the following conjectural properties.
\begin{enumerate}

\item $\T(X)$ is a (neutral) tannakian category over $\Q$ with a fiber functor
$\w:\T(X) \to \Vec_\Q$ to the category of finite dimensional rational vector
spaces.

\item Each object $M$ of $\T(X)$ has an increasing filtration called the weight
filtration
$$
\cdots \subseteq W_{m-1}M \subseteq W_m M \subseteq W_{m+1}M \subseteq \cdots,
$$
whose intersection is $0$ and whose union is $M$. Morphisms of $\T(X)$ should
be strictly compatible with the weight filtration --- that is,  the functor
$$
\Gr^W_\dot:M \mapsto \bigoplus_m \Gr^W_m M :=
\bigoplus_m W_mM/W_{m-1}M
$$ 
to graded objects in $\T(X)$ should be an exact tensor functor.

\item $\T(X)$ contains ``the Tate motive $\Q(1)$'' over $\Spec R$ where $R$ is
the base ring (here either a field or $\OFS$). This can be considered as the
dual of the local system  $R^1 f_*(\Q)$ over $\Spec R$, where  $f$ is the
structure morphism of the multiplicative group $\Gm$, i.e., $f:\Gm \otimes R
\to \Spec R$. Put $\Q(n) := \Q(1)^{\otimes n}$ for $n\in \Z$. (Negative tensor
powers are defined by duality.)

\item There should be {\it realization functors} to various categories such as
$\l$-adic \'etale local systems over $X[1/\l]:= X\otimes_R R[1/\l]$ (where $R$
is the base ring and $\l$ does not divide the characteristic of $R$),
variations of mixed Hodge structure over $X$, etc. These functors should be
faithful, exact tensor functors. These functors are related by natural
comparison transformations. The Betti, de~Rham, $\l$-adic and crystalline
realizations of $\Q(1)$ should be the Betti, de~Rham, $\l$-adic and crystalline
versions of $H_1(\Gm)$.

\item For each object $M$, $\Gr^W_{2m+1} M$ is trivial and $\Gr^W_{2m} M$ is
isomorphic to the direct sum of a finite number of copies of $\Q(-m)$.

\end{enumerate}
The last property characterizes {\em mixed Tate} motives among mixed motives.
The category $\T(X)$, being tannakian, is equivalent to the category of finite
dimensional representations of a proalgebraic $\Q$-group  $\pi_1(\T(X),\w)$,
which represents the tensor automorphism group of the fiber functor $\w$. We
denote it simply by $\pi_1(\T(X))$, if  the selection of $\w$ does not matter.

There are several approaches to constructing the category $\T(X)$, at least
when $X$ is the spectrum of a field or $X=\XFS$, such as those of Bloch-Kriz
\cite{bloch-kriz}, Levine \cite{levine:tate}, and
Goncharov\cite{goncharov:mtm}.

We follow Deligne \cite{deligne} and Jannsen \cite{jannsen}, who define a
motive to be a ``compatible set of realizations'' of ``geometric origin.'' This
is a tannakian category. Deligne does not define what it means to be of
geometric origin, but wants it to be broad enough to include those compatible
realizations that occur in the unipotent completion of fundamental groups of
varieties in addition to subquotients of cohomology groups. We refer the reader
to Section~1 of Deligne's paper \cite{deligne} for the definition of compatible
set of realizations. One example is $\Q(1)$, defined as $H_1(\Gm_{/\Z})$,
another is the extension $E_x$ of $\Q(0)$ by $\Q(n)$ coming from $x\in
K_{2n-1}(\XFS)$ described in the previous section.

The hope is that
$$
\Hmot^m(X,\Q(n)) \cong \Ext^m_{\T(X)}(\Q(0),\Q(n))
$$
holds when $X$ is the spectrum of a field or $\XFS$.%
\footnote{
If this is true, then $\Hmot^m(\Spec F,\Q(n))$ will vanish when $n < 0$ and
$m=0$, and when $n \le 0$ and $m>0$. This vanishing is a conjecture of
Beilinson and Soul\'e. It is known for number fields, for example.
} 
This covers the cases of interest for us.  In general, one expects that
motivic cohomology groups of $X$ can be computed as hyper-exts:
$$
\Hmot^m(X,\Q(n)) \cong \bH^m(X,\sExt_\T^\dot(\Q(0),\Q(n))).
$$

Deligne's conjecture (Conjecture~\ref{del_conj})
will be a consequence of:

\begin{postulate}
\label{post}
If $X=\XFS$, there is a category of mixed Tate motives $\T(X)$ over $X$ with
the above mentioned properties. It has the property that there is a natural
isomorphism
$$
\Hmot^m(X,\Q(n)) \cong \Ext^m_{\T(X)}(\Q(0),\Q(n)),
$$
which is compatible with Chern maps.
\end{postulate}
\subsection*{Examples of Mixed Tate Motives over $\Spec \Z$}

One of the main points of \cite{deligne} is to show that the unipotent
completion of the fundamental group of $\Pminus$ is an example of a
mixed Tate motive (actually a pro-mixed Tate motive), smooth over
$\Spec \Z$.

As base point, take $\v$, the tangent vector of $\P^1$ based at $0$ that
corresponds to $\partial/\partial t$, where $t$ is the natural parameter on
$\Pminus$. Deligne \cite{deligne} shows that the unipotent completion
of $\pi_1(\Pminus,\v)$ is a mixed Tate motive over $\Spec \Z$ by exhibiting
compatible Betti, \'etale, de~Rham and crystalline realizations of it. It is
smooth over $\Spec \Z$ essentially because the pair $(\Pminus, \v)$ has
everywhere good reduction.

There is an interesting relation to classical polylogarithms which was
discovered by Deligne (cf.\ \cite{deligne}, \cite{beilinson-deligne},
\cite{hain:polylogs}). There is a polylog local system $P$, which is a motivic
local system over $\P^1_{\Z}-\{0,1\infty\}$ in the point of view of
compatible realizations. Its Hodge-de~Rham realization is a variation of mixed
Hodge structure over the complex points of $\Pminus$ whose periods are given by
$\log x$ and the classical polylogarithms: $\Li_1(x) = - \log(1-x)$, $\Li_2(x)$
(Euler's dilogarithm), $\Li_3(x)$, and so on. Here $\Li_n(x)$ is the
multivalued holomorphic function on $\PminusC$ whose principal branch is given
by
$$
\Li_n(x) = \sum_{k=1}^\infty \frac{x^k}{k^n}
$$
in the unit disk.\footnote{This goes back to various letters of Deligne.
Accounts can be found, for example, in \cite{beilinson-deligne} and
\cite{hain:polylogs}.}

The fiber $P_\v$ of $P$ over the base point $\v$ is a mixed Tate motive over
$\Spec \Z$ and has periods the values of the Riemann zeta function at
integers $n>1$. In fact, $P_\v$ is an extension
$$
0 \to \bigoplus_{n \ge 1} \Q(n) \to P_\v \to \Q(0) \to 0
$$
and thus determines a class
$$
(e_n)_n \in \bigoplus_{n\ge 1} \Ext^1_\H(\Q(0),\Q(n)).
$$
The class $e_n$ is trivial when $n=1$ and is the coset of $\zeta(n)$
in
$$
\C/(2\pi i)^n \Q \cong \Ext^1_\H(\Q(0),\Q(n))
$$
when $n>1$. Since $\zeta(2n)$ is a rational multiple of $\pi^{2n}$, each
$e_{2n}$ is trivial.

Deligne computes the $\l$-adic realization of $P_\v$ in \cite{deligne}
and shows that the polylogarithm motive is a canonical quotient of the
enveloping algebra of the Lie algebra of the unipotent completion of
$\pi_1(\Pminus,\v)$. (See also \cite{beilinson-deligne} and
\cite{hain:polylogs}.)

\section{The Motivic Lie Algebra of $\XFS$ and Deligne's Conjectures}

Assume that $X$ is as in Section~\ref{sec:mtm}, and that there is a category of
mixed Tate motives $\T(X)$ with properties (i)--(v) in Section~\ref{sec:mtm}.
Since $\T(X)$ is tannakian, it is determined by its tannakian fundamental group
$\pi_1(\T(X))$, which is  an extension of $\Gm$ by a prounipotent $\Q$-group
\begin{equation}\label{eq:extension}
1 \to \U_X \to \pi_1(\T(X)) \to \Gm \to 1
\end{equation}
as we shall now explain.

The category of pure Tate motives is the tannakian subcategory of $\T(X)$
generated by $\Q(1)$. By the faithfulness of realization functors, it is
equivalent to the category of finite dimensional graded $\Q$-vector spaces, 
and hence to the category of finite dimensional representations of $\Gm$;
$\Q(n)$ corresponds to the $n$th power of the  standard representation. This
induces a group homomorphism between the tannakian  fundamental groups
$$
\pi_1(\T(X)) \to \pi_1(\text{pure Tate motives})\cong \Gm.
$$
Since the category of pure Tate motives is a full subcategory and every
subobject of a pure Tate motive is pure, this morphism is surjective (cf.\
\cite[Proposition~2.21a]{dmos}), and the properties of the weight filtration 
imply the unipotence of its kernel $\U_X$, thus we have (\ref{eq:extension}).
The Lie algebra $\t_X$ of $\pi_1(\T(X))$ is an extension
$$
0 \to \u_X \to \t_X \to \Q \to 0
$$
where $\u_X$ is pronilpotent. This $\t_X$ is called the motivic Lie algebra of
$X$. We shall see that the knowledge of the cohomologies of $\u_X$ (as
$\Gm$-modules) is equivalent to the knowledge of the extension  groups
$\Ext^\dot_{\T(X)}(\Q(0),\Q(m))$ for all $m$ in the next section.

\subsection{Extension groups in a tannakian category}
We start with a general setting.  Let $K$ be a field of characteristic zero.
Let $\cG$ be a proalgebraic group (in this paper a proalgebraic group means an
affine proalgebraic group) over $K$, or equivalently, an affine group scheme
over $K$ (cf.\ \cite{dmos}). A $\cG$-module $V$ is a (possibly infinite
dimensional) $K$-vector space with algebraic $\cG$-action (cf.\
\cite{Jantzen}). The category of $\cG$-modules is abelian with enough
injectives,  and hence we have the cohomology groups
$$
H^m(\cG,V):=\Ext^m_\cG(K,V)
$$
defined as the extension groups, where $K$ denotes the trivial representation.
The right hand side has an interpretation as Yoneda's extension groups, i.e.,
as the set of equivalence classes of $m$-step extensions (see \cite{yoneda2}).
Since each $\cG$-module is locally finite \cite[2.13]{Jantzen}, every $m$-step
extension representing an element of $\Ext^m_\cG(K,V)$ can be replaced by an
equivalent extension consisting of finite dimensional modules when $V$ is
finite dimensional. Thus,  the right hand side does not change when the
category of $\cG$-modules is replaced by the category of finite dimensional
$\cG$-modules.

Let $\T$ be a neutral tannakian category over $K$ with a fiber functor $\w$,
and let $\cG$ be its tannakian fundamental group with base point $\w$.  Since
$\T$ is isomorphic to the category of finite dimensional $\cG$-modules, we have
the following.

\begin{lemma}
Let $\T$ be a neutral tannakian category and $\cG$ be its tannakian fundamental
group. Then, for any object $V$, we have
$$
\Ext^m_\T(K,V)\cong H^m(\cG,V). \qed
$$
\end{lemma}

Suppose that $\cG$ is an extension
$$
1 \to \U \to \cG \to R \to 1
$$
of proalgebraic groups over $K$. Then, for any $\cG$-module $V$, we have  the
Lyndon-Hochschild-Serre spectral sequence  (cf.\ \cite[6.6
Proposition]{Jantzen}):
$$
E_2^{s,t}=H^s(R,H^t(\U,V)) \Rightarrow H^{s+t}(\cG,V).
$$
If $R$ is a reductive algebraic group,  then every $R$-module is completely
reducible. Consequently, $H^s(R,V)$ vanishes for $s\geq 1$ for all $V$, and
$$
H^m(\cG,V)\cong H^0(R,H^m(\U,V)).
$$
If, in addition, the action of $\cG$ on $V$ factors through $R$, then one has
an $R$-module isomorphism
$$
H^m(\U,V)\cong H^m(\U,K)\otimes V.
$$

Moreover, if we assume that $\U$ is prounipotent, then its Lie algebra $\u$ is
a projective limit
$$
\u\cong \limproj{\n} \u/\n
$$
of finite dimensional nilpotent Lie algebras. It has a topology as a projective
limit, where each $\u/\n$ is viewed as a  discrete topological space.

Let $V$ be a continuous $\u$-module over $K$. The continuous cohomology
$\Hcts^m(\u,V)$ is defined as  the extension group $\Ext^m(K,V)$ in the
category  of continuous $\u$-modules. We denote  $\Hcts^m(\u,K)$ by
$\Hcts^m(\u)$. It is easy to show that
$$
\Hcts^m(\u) \cong \liminj{\n}H^m(\u/\n),
$$
where $H^m(\u/\n)$ can be computed as the cohomology of the complex of cochains
$$
\Hom(\Lambda^\dot (\u/\n), K).
$$
The following is standard.

\begin{proposition}
Let $\u$ be a pronilpotent Lie algebra, and let $H_1(\u)$ denote the
abelianization of $\u$. Then
$$
H_1(\u) \cong \Hom(\Hcts^1(\u),K).
$$
If $H^2(\u)=0$, then $\u$ is free.
\end{proposition}

It is also well known that the category of $\U$-modules is equivalent to the
category of continuous $\u$-modules. Hence we have
$$
H^m(\U,K)\cong \Hcts^m(\u).
$$
Putting this together, we have the following.

\begin{theorem}\label{Th:unip-cohomology}
Suppose that $1 \to \U \to \cG \to R \to 1$ is a short exact sequence of
pro-algebraic groups over a field $K$ of  characteristic zero. Assume that $R$
is a reductive  algebraic group, and that $\U$ is a prounipotent group. Let
$\u$ be the Lie algebra of $\U$. If $V$ is an $R$-module, considered as a
$\cG$-module, then
$$
H^m(\cG,V) \cong (\Hcts^m(\u)\otimes V)^R.
$$
Consequently, we have the $R$-module isomorphism
$$
\Hcts^m(\u) \cong \bigoplus_\alpha (H^m(\cG,V_\alpha) \otimes V_\alpha^*),
$$
where $\{V_\alpha\}$ is a set of representatives of the isomorphism classes of
irreducible $R$-modules, and $(\blank)^\ast$ denotes $\Hom(\blank,K)$.
\end{theorem}

\subsection{Deligne's conjecture}
\label{discussion}
By applying Theorem~\ref{Th:unip-cohomology} to (\ref{eq:extension}), we have
$$
\Ext^m_{\T(X)}(\Q(0),\Q(n)) \cong [\Hcts^m(\u_X)\otimes\Q(n)]^\Gm
$$
and $\Gm$-module isomorphisms
\begin{equation}\label{eq:ux}
\Hcts^m(\u_X)\cong
\bigoplus_{n \in \Z}\Ext^m_{\T(X)}(\Q(0),\Q(n))\otimes \Q(-n),
\end{equation}
where $\u_X$ is the Lie algebra of $\U_X$, the prounipotent radical of
$\pi_1(\T(X),\w)$. By a weight argument, each extension on the right hand side
vanishes if $n\leq m-1$. Postulate~\ref{post} says that these $\Ext$ groups
should be the motivic cohomology groups of $X$, and Theorem~\ref{mot} says that
they should be isomorphic to the Adams eigenspaces of the $K$-groups of $X$:

\begin{proposition}
\label{mtm}
\label{invariants}
Assume the existence of a category $\T(\XFS)$ of mixed Tate motives over $\XFS$
with properties (i)--(v) as in Section~\ref{sec:mtm}. Suppose that
Postulate~\ref{post} holds for all $n\geq 1$. Let $\U_\XFS$ be the unipotent
radical of  $\pi_1(\T(\XFS),\w)$, and $\u_\XFS$ be its Lie algebra. Then there
is a natural $\Gm$-module isomorphism
$$
\Hcts^1(\u_\XFS)\cong
\bigoplus_{n\geq 1} K_{2n-1}(\XFS)\otimes_\Z \Q(-n),
$$
and $\Hcts^m(\u_\XFS)=0$ whenever $m\geq 2$.
Moreover, the exactness of $\Gr^W_\dot$ implies that
$$
H_1(\Gr^W_\dot\u_\XFS) = \bigoplus_{n\geq 1} K_{2n-1}(\XFS)^\ast \otimes \Q(n)
$$
and that
$$
H^m(\Gr^W_\dot\u_\XFS) = 0 \text{ when } m>1.
$$
It follows from this that $\Gr^W_\dot\u_\XFS$ is isomorphic to the free
Lie algebra generated by $H_1(\Gr^W_\dot\u_\XFS)$.
\end{proposition}

Let us assume that there is a category $\T(\XFS)$ satisfying (i)--(v)
in Section~\ref{sec:mtm}. Then Postulate~\ref{post} is equivalent to
the following conjecture of Deligne:

\begin{conjecture}[Deligne]
\label{del_conj}
\begin{enumerate}
\item{\cite[8.2.1]{deligne}}
For the category $\T(\XFS)$ of motives smooth over $\XFS$ one has a natural
isomorphism
$$
\Ext_{\T(\XFS)}^1(\Q(0),\Q(n)) \cong 
K_{2n-1}(\XFS)\otimes \Q \text{ for all }n,
$$
which is compatible with the Chern mappings.
\item{\cite[8.9.5]{deligne}} The group $\pi_1(\T(\XFS))$ is an extension
of $\Gm$ by a free prounipotent group.
\end{enumerate}
\end{conjecture}

Note that by Definition~\ref{def}, Theorem~\ref{mot} 
and the isomorphism (\ref{eq:ux}),
(i) is equivalent to Postulate~\ref{post} for $m=1$, and that  
(ii) is equivalent to Postulate~\ref{post} for $m\geq 2$.

\subsection*{Consequences of Deligne's Conjecture}

Deligne's conjecture suggests restrictions on the action of Galois groups
on pro-$\l$ completions of fundamental groups of curves. Here is a sketch of
how this should work. 

As in the beginning of Section~\ref{sec:mtm}, there should be a Betti
realization functor
$$
\real_B : \T(\XFS) \to \{\text{$\Q$-vector spaces}\}
$$
to the category of $\Q$-vector spaces, and an $\l$-adic realization functor
$$
\real_\l : \T(\XFS) \to \{\text{$\l$-adic $G_F$-modules}\},
$$
to the category of the $\Ql$-vector spaces with a continuous $G_F$-action.  The
Galois modules should be unramified outside $S\cup [\l]$, where $[\l]$ denotes
the set of primes of $F$ over $\l$. We choose $\real_B$ as our fiber functor
$\w$. Let $\w_\l$ denote the functor $\real_\l$ which forgets the $G_F$-action.
Conjecturally, there is a comparison isomorphism
$$
\w \otimes \Ql \cong \w_\l,
$$
so we shall identify these two. Define $\T(\XFS)\otimes \Ql$ to be the
tannakian category whose objects are the same as those of $\T(\XFS)$ and whose
hom-sets are those of $\T(\XFS)$ tensored with $\Ql$. The $\l$-adic
realization functor induces a functor
\begin{equation}\label{realization}
\real_\l : \T(\XFS)\otimes \Ql \to \{\text{$\l$-adic $G_F$-modules}\}
\end{equation}
(by an abuse of notation we denote it by $\real_\l$ again), and by forgetting
the Galois action a fiber functor $\w_\l:\T(\XFS)\otimes \Ql\to \Vec_\Ql$
(under a similar abuse of notation). Through the comparison isomorphism, it is
easy to show that
$$
\pi_1(\T(\XFS)\otimes\Ql,\w_\l)\cong \pi_1(\T(\XFS),\w)\otimes\Ql.
$$
The following is closely related to the Tate conjecture on Galois modules.
\begin{postulate}\label{post2}
The realization functor $\real_\l$ in (\ref{realization}) is 
fully faithful, and its image is closed under taking subobjects.
\end{postulate}

The first condition is that every Galois compatible morphism comes from a
morphism of motives up to extension of scalars, and the second condition is
that every Galois submodule arises as an $\l$-adic
realization of a motive. We shall see that this
postulate follows from Deligne's Conjecture~\ref{del_conj} and our
Theorem~\ref{main} (see Corollary~\ref{cor:post2-del}).

Every element of $G_F$ gives an automorphism of the 
forgetful fiber functor of the category of $G_F$-modules
(i.e.\ forgetting the Galois action), and hence an automorphism of $\w_\l$.
Thus we have a homomorphism
\begin{equation}\label{GFStotannaka}
G_F \to 
\pi_1(\T(\XFS)\otimes\Ql,\w_\l)(\Ql) \cong \pi_1(\T(\XFS),\w)(\Ql).
\end{equation}
In addition, the $G_F$-action on  the $\l$-adic realization of any
(pro)object $M$ of $\T(\XFS)$ factors through
$\pi_1(\T(\XFS))\otimes\Ql$ via the morphism (\ref{GFStotannaka}).

\begin{proposition}\label{prop:zar}
Postulate~\ref{post2}  is equivalent to the statement that the above morphism
(\ref{GFStotannaka}) has Zariski dense image.
\end{proposition}

\begin{proof}
Let $\cG$ denote the tannakian fundamental group of the category of finite
dimensional $\l$-adic $G_F$-modules. By \cite[Prop.~2.21a]{dmos}, the
conditions in  Postulate~\ref{post2} are equivalent to the surjectivity of $\cG
\to \pi_1(\T(\XFS)\otimes\Ql,\w_\l)$. It is  a general fact that the image of
$G_F \to \cG(\Ql)$ is Zariski dense.
\end{proof}

Assuming Postulate~\ref{post2}, the Zariski density of the image of
(\ref{GFStotannaka}) implies that for any object $M$ of $\T(\XFS)$, the Zariski
closure of the image of $G_F$ in $\Aut(M)$  should be a quotient of
$\pi_1(\T(\XFS))\otimes\Ql$. We can define a filtration $J_M^\dot$ on $G_F$
(which depends on $M$) by
$$
J_M^m G_F := \text{ the inverse image of } W_m\Aut M.
$$
The image of the Galois group $G_{F(\mu_{\l^\infty})}$ of
$F(\mu_{\l^\infty})$ in $\pi_1(\T(\XFS),\w_\l)$ will lie in its prounipotent
radical and should be Zariski dense in it. The exactness of $\Gr^W_\dot$
will then imply that
$$
\big(\bigoplus_{m<0}\Gr_{J_M}^m G_F\big)\otimes_\Zl \Ql
$$
(a Lie algebra) is a quotient of $\Gr^W_\dot \u_\XFS$, and hence generated
by
$$
\bigoplus_{m\ge 1} \Hom(K_{2m-1}(\XFS), \Ql(m)).
$$

For example, the pronilpotent Lie algebra $\p$ of the unipotent completion of
$\pi_1(\Pminus,\v)$ should be a pro-object of $\T(\Spec \Z)$. One should
therefore expect that the graded Lie algebra
$$
\big(\bigoplus_{m<0}\Gr_{J_{\p}}^m G_\Q\big)\otimes_\Zl \Ql
$$
is generated by elements $z_3,z_5,z_7,\dots$, where $z_m$ has weight $-2m$.

Following Ihara \cite{Ihara1}, we define
$$
I_\l^m G_\Q = \ker\{G_\Q \to \Out(\pi_1^\prol(\PminusC,\v)/L^{m+1})\}
$$
where $L^m$ denotes the $m$th term of the lower central series of the pro-$\l$
completion of $\pi_1(\PminusC,\v)$. This is related to the filtration
$J_{\p}^\dot$ by
$$
I_\l^m G_\Q = J_{\p}^{-2m} G_\Q = J_{\p}^{-2m+1} G_\Q.
$$
Making this substitution, we are led to the following conjecture, stated
by Ihara in \cite[p.~300]{Ihara3} and which he attributes to Deligne.

\begin{conjecture}[Deligne]
\label{del-ihara_conj}
The Lie algebra
$$
\big[\bigoplus_{m>0} \Gr^m_\I G_\Q\big]\otimes \Ql
$$
is generated by generators $s_3, s_5, s_7, \dots$, where
$s_m \in \Gr^m_\I G_\Q$.
\end{conjecture}

Deligne also asked whether this Lie algebra is free. A related conjecture of
Goncharov \cite[Conj.~2.1]{goncharov}, stated below, and the questions of
Drinfeld \cite{drinfeld} can be `derived' from Deligne's
Conjecture~\ref{del_conj} in a similar fashion. The freeness questions are more
optimistic and are equivalent to the statement that the representation of the
motivic Galois group $\pi_1(\T(\Spec \Z))$ in the automorphisms of the
$\l$-adic unipotent completion of the fundamental group of $\Pminus$ is
faithful.  The computational results \cite{Ihara3}, \cite{matsumoto} and
\cite{tsunogai} give support to the belief that this Lie algebra is free.
Indeed, these computations show that $\Gr_{\I}^{>0} G_\Q$ is free up to
$\Gr_\I^{12}G_\Q$.

\begin{conjecture}[Goncharov]
\label{gonch_conj}
The Lie algebra of the Zariski closure of the Galois group of
$\Q(\mu_{\l^\infty})$ in the automorphism group of the $\l$-adic
unipotent completion of $\pi_1(\Pminus,\v)$ is a prounipotent Lie algebra
freely generated by elements $z_3, z_5, z_7, \dots$, where $z_m$ has weight
$-2m$.
\end{conjecture}

Deligne's Conjecture~\ref{del-ihara_conj} above and the generation part of
Goncharov's conjecture are proved in \cite{hain-matsumoto}. A brief sketch
of their proofs is  given in Section~\ref{proof}. Modulo technical details,
the main point is the computation of the tannakian fundamental group of
the candidate for $\T(\XFS)\otimes\Ql$ given in the next section.

\subsection*{Polylogarithms Revisited}
Assuming the existence of $\t_{\Spec \Z}$
(i.e.\ the Lie algebra of $\pi_1(\T(\Spec \Z))$), 
we can give another interpretation
of the fiber $P_\v$ of the polylogarithm local system. Being a motive
over $\Spec \Z$, it is a $\t_{\Spec \Z}$-module. Note that since
$$
W_{-1} P_\v = \bigoplus_{n\ge 1} \Q(n),
$$
a direct sum of Tate motives (no non-trivial extensions), the restriction
of the $\t_{\Spec\Z}$-action on $ W_{-1} P_\v$ to $\u_{\Spec \Z}$ is trivial
and $[\u_{\Spec \Z},\u_{\Spec \Z}]$ annihilates $P_\v$.
Since $P_\v$ is an extension of $\Q(0)$ by $W_{-1} P_\v$, this implies that
there is a homomorphism
$$
\psi : \big(\t_{\Spec \Z}/[\u_{\Spec \Z},\u_{\Spec \Z}]\big)\otimes \Q(0)
\longrightarrow P_\v
$$
such that the diagram
$$
\xymatrix{
\t_{\Spec \Z} \otimes P_\v \ar[d]_{\text{quotient}}
\ar[r]^{\quad \text{action}} & P_\v \cr
\t_{\Spec \Z}/[\u_{\Spec \Z},\u_{\Spec \Z}] \ar[ur]_\psi}
$$
commutes. By comparing graded quotients, it follows that $\psi$ is an
isomorphism
$$
P_\v \cong \t_{\Spec \Z}/[\u_{\Spec \Z},\u_{\Spec \Z}]
$$
of motives over $\Spec \Z$.

\section{$\l$-adic Mixed Tate Modules over $\XFS$}
\label{wtd_ell_mods}

In this section, we describe a candidate for the category of $\l$-adic
realizations of objects and morphisms of $\T(\XFS)$. This is essentially the
category constructed by Deligne and Beilinson in their unpublished manuscript
\cite{beilinson-deligne:orig}. It is purely Galois-representation theoretic,
and requires no postulates.  For technical reasons, we assume that $S$
contains  $[\l]$, the set of primes over $\l$. This condition will be removed
in Section~\ref{sec:crys}. By a finite dimensional $\GFS$-module, we shall mean
a finite dimensional $\Ql$-vector space on which $\GFS$ acts continuously.

We define the category $\T_\l(\XFS)$ of $\l$-adic mixed Tate modules which are
smooth over $\XFS$ to be the category whose objects are finite dimensional
$\GFS$-modules $M$ that are equipped with a {\it weight filtration}
$$
\cdots \subseteq W_{m-1}M \subseteq W_m M \subseteq W_{m+1}M \subseteq \cdots
$$
of $M$ by $\GFS$-submodules. The weight filtration satisfies:
\begin{enumerate}
\item all odd weight graded quotients of $M$ vanish: $\Gr^W_{2m+1} M = 0$;
\item $\GFS$ acts on its $2m$th graded quotient $\Gr^W_{2m} M$ via
the $(-m)$th power of the cyclotomic character,
\item the intersection of the $W_m M$ is trivial and their union is all
of $M$.
\end{enumerate}
Morphisms are $\Ql$-linear, $\GFS$-equivariant mappings. These will
necessarily preserve the weight filtration, so that $\T_\l(\XFS)$ is a full
subcategory of the category of $\GFS$-modules.

The category $\T_\l(\XFS)$ is a tannakian category over $\Ql$ with a fiber
functor $\w'$ that takes an object to its underlying $\Ql$-vector space. We
shall denote the tannakian fundamental group of this category by
$\AFSl:=\pi_1(\T_\l(\XFS),\w')$. 
Every element of $\GFS$ acts on $\w'$, 
which induces a natural, continuous homomorphism
$$
\rho : \GFS \to \AFSl(\Ql).
$$
This has Zariski-dense image as $\T_\l(\XFS)$ is a full subcategory of the
category of $\GFS$-modules, closed under taking subobjects 
(cf.\ \cite[Proposition~2.21a]{dmos}).

\subsection*{Relation to Mixed Tate Motives over $\XFS$}
As explained in Section~\ref{sec:mtm}, the existence of a category $\T(\XFS)$
of mixed Tate motives over $\XFS$ satisfying (i)--(v) in Section~\ref{sec:mtm}
implies the existence of an $\l$-adic realization functor
$$
\real_\l : \T(\XFS)\otimes \Ql \to \T_\l(\XFS).
$$
This will induce a morphism of tannakian fundamental groups
$$
\AFSl=\pi_1(\T_\l(\XFS),\w') \to \pi_1(\T(\XFS),\w)\otimes \Ql.
$$

The main result of \cite{hain-matsumoto} may be interpreted as saying that
$\AFSl=\pi_1(\T_\l(\XFS),\w')$ is isomorphic to the conjectured value of the
$\Ql$-form $\pi_1(\T(\XFS),\w)\otimes\Ql$ of the motivic fundamental
group of $\XFS$. We shall explain this in Section~\ref{compn}.

It is interesting to note that we have not restricted to $\GFS$-modules of
geometric origin as Deligne would like to. So one consequence of our result
is that, if Deligne's Conjecture~\ref{del_conj} is true, then all weighted
$\l$-adic $\GFS$-modules and their morphisms will be of geometric origin.

\section{Weighted Completion of Profinite Groups}

In this and the subsequent two sections we will sketch how to compute the
tannakian fundamental group $\pi_1(\T_\l(\XFS),\w')$ of the category of
$\l$-adic mixed Tate modules smooth over $\XFS$, which was defined in
Section~\ref{wtd_ell_mods}. It is convenient to work in greater
generality.\footnote{ We may generalize further:  weighted completion and its
properties in  this section are unchanged even if we replace $\Ql$ by an
arbitrary topological field of characteristic zero and $\G$ by an arbitrary
topological group.}

Suppose that $R$ is a reductive algebraic group  over $\Ql$ and that $w : \Gm
\to R$ is a central cocharacter --- that is, its image is contained in the
center of $R$. It is best to imagine that $w$ is non-trivial as the theory of
weighted completion is uninteresting if $w$ is trivial.

Suppose that $\G$ is a profinite group and that a homomorphism $\rho : \G \to
R(\Ql)$ has Zariski dense image and  is continuous where we view $R(\Ql)$ as an
$\l$-adic Lie group.

By a {\it weighted} $\G$-module with respect to $\rho$ and $w$ we shall mean a
finite dimensional $\Ql$-vector space with continuous $\G$-action together with
a weight filtration
$$
\cdots \subseteq W_{m-1}M \subseteq W_m M \subseteq W_{m+1}M \subseteq \cdots
$$
by $\G$-invariant subspaces. These should satisfy:
\begin{enumerate}
\item the intersection of the $W_mM$ is $0$ and their union is $M$,
\item for each $m$, the representation $\G \to \Aut \Gr^W_m M$ should
factor through $\rho$ and a homomorphism $\phi_m : R \to \Aut \Gr^W_m M$,
\item $\Gr^W_m M$ has weight $m$ when viewed as a $\Gm$-module via
$$
\begin{CD}
\Gm @>w>> R @>{\phi_m}>> \Aut \Gr^W_m M.
\end{CD}
$$
That is, $\Gm$ acts on $\Gr^W_m M$ via the $m$th power of the standard
character.
\end{enumerate}
The category of weighted $\G$-modules consists of the $\G$-equivariant
morphisms between weighted $\G$-modules. These morphisms automatically
preserve the weight filtration and are strict with respect to it; that is,
the functor $\Gr^W_\dot$ is exact.

One can show that the category of weighted $\G$-modules is tannakian,
with fiber functor $\w'$ given by forgetting the $\G$-action.

\begin{definition}
The weighted completion of $\G$ with respect to $\rho : \G \to R(\Ql)$
and $w : \Gm \to R$ is the tannakian fundamental group of the category of
weighted $\G$ modules with respect to $\rho$ and $w$.
\end{definition}

Denote the weighted completion of $\G$ with respect to $\rho$ and $w$ by
$\cG$. There is a natural homomorphism $\G \to \cG(\Ql)$ which has Zariski
dense image as we shall see below.

This definition differs from the one given in \cite[Section~5]{hain-matsumoto},
but is easily seen to be equivalent to it. (See below.) In particular, we can
apply it when:
\begin{itemize}
\item $\G$ is $\GFS$,
\item $R$ is $\Gm$ and $w : \Gm \to \Gm$ takes $x$ to $x^{-2}$,
\item $\rho$ is the composite of the $\l$-adic cyclotomic character
$\chi_\l : \GFS \to \Zlx$ with the inclusion $\Zlx \hookrightarrow \Qlx$.
\end{itemize}
In this case, the category of weighted $\G$-modules
is nothing but the category of mixed Tate modules $\T_\l(\XFS)$.
Recall that we denote the corresponding weighted completion by 
$\AFSl:=\pi_1(\T_\l(\XFS),\w')$.

\subsection*{Equivalence of Definitions}
Here we show that the definition of weighted completion given in
\cite{hain-matsumoto} agrees with the one given here.

Suppose that $G$ is a linear algebraic group over $\Ql$ which is an extension
$$
1 \to U \to G \to R \to 1
$$
of $R$ by a unipotent group $U$. Note that $H_1(U)$ is naturally an $R$-module,
and therefore a $\Gm$-module via the given central cocharacter $w:\Gm \to R$.
We can decompose $H_1(U)$ as a $\Gm$-module:
$$
H_1(U) = \bigoplus_{n\in\Z} H_1(U)_n
$$
where $\Gm$ acts on $H_1(U)_n$ via the $n$th power of the standard character.
We say that $G$ is a {\it negatively weighted extension of $R$}
if $H_1(U)_n$ vanishes whenever $n\ge 0$.

Given a continuous homomorphism $\rho : \G \to R(\Ql)$ with Zariski dense
image, we can form a category of pairs $(\rhotilde, G)$, where $G$ is a
negatively weighted extension of $R$ and $\rhotilde : \G \to G(\Ql)$ is
a continuous homomorphism 
that lifts $\rho$. Morphisms in this category
are given by homomorphisms between the $G$s that respect the projection
to $R$ and the lifts $\rhotilde$ of $\rho$. The objects of this category,
where $\rhotilde$ is Zariski dense, form an inverse system. Their inverse
limit is an extension
$$
1 \to \U \to \cG \to R \to 1
$$
of $R$ by a prounipotent group. There is a natural homomorphism $\rhohat : \G
\to \cG(\Ql)$, which is continuous in a natural sense. It has the following
universal mapping property: if $\rhotilde : \G \to G(\Ql)$ is an object of this
category, then there is a unique homomorphism $\phi : \cG \to G$ that commutes
with the projections to $R$ and with the homomorphisms $\rhotilde : \G \to
G(\Ql)$ and $\rhohat : \G \to \cG(\Ql)$. In \cite{hain-matsumoto}, the weighted
completion is defined to be this inverse limit. The equivalence of the two
definitions follows from the following result.

\begin{proposition}
The inverse limit above is naturally isomorphic to the weighted completion of
$\G$ relative to $\rho$ and $w$.
\end{proposition}

\begin{proof}
Denote the inverse limit by $\cG$ and by $\M=\M(\rho,w)$ the category of
weighted $\G$-modules with respect to $\rho : \G \to R(\Ql)$ and $w$. We will
show that $\M$ is the category of finite dimensional $\cG$-modules, from which
the result follows.

Suppose that $M$ is an object of $\M$. Then the Zariski closure of $\G$
in $\Aut M$ is an extension
$$
1 \to U \to G \to R' \to 1
$$
of a quotient of $R$ by a unipotent group. Here
$R'$ is the Zariski closure of the image of $\G$
in $\Aut \Gr^W_\dot M$. Because the action of $\G$ on
each weight graded quotient factors through $\rho$, and because $\Gm$ acts
on the $m$th weight graded quotient of $M$ with weight $m$, it follows that
this is a negatively weighted extension of $R'$. Pulling back this extension
along the projection $R \to R'$, we obtain a negatively weighted extension
$$
1 \to U \to \widetilde{G} \to R \to 1
$$
of $R$ and a continuous homomorphism $\G \to \widetilde{G}(\Ql)$ that lifts
both $\rho$ and the homomorphism $\G \to R'(\Ql)$. By the universal mapping
property of $\cG$, there is a natural homomorphism $\cG \to \widetilde{G}$
compatible with the projections to $R$ and the homomorphisms from $\G$
to $\cG(\Ql)$ and $\widetilde{G}(\Ql)$. Thus every object of $\M$ is naturally
a $\cG$-module. It is also easy to see that every morphism of $\M$ is
$\cG$-equivariant.

Conversely, suppose that $M$ is a finite dimensional $\cG$-module. Composing
with the natural homomorphism $\rhohat : \G \to \cG(\Ql)$ gives $M$ the
structure of a $\G$-module. In \cite[Sect.~4]{hain-matsumoto}, it is proven
that every $\cG$-module has a natural weight filtration with the property
that the action of $\cG$ on each weight graded quotient factors through the
projection $\cG \to R$ and that $\Gm$ acts with weight $m$ on the $m$th
weight graded quotient. It follows that $M$ is naturally an object of
$\M$. Since $\cG$-equivariant mappings are naturally $\G$-equivariant,
this proves that $\M$ is naturally the category of finite dimensional
$\cG$-modules, which completes the proof.
\end{proof}

\section{Computation of Weighted Completions}

Suppose that $\G$, $R$, $\rho : \G \to R(\Ql)$ and $w : \Gm \to R$ are as
above.  The weighted completion $\cG$ of $\G$ is controlled by the
low-dimensional cohomology groups $\Hcts^\dot(\G,V)$ of $\G$ with coefficients
in certain irreducible representations $V$ of $R$. If one knows these
cohomology groups, as we do in the case of $\GFS$, one can sometimes determine
the structure of the weighted completion. These cohomological results are
stated in this section.

The weighted completion of $\G$ with respect to $\rho$ and $w$ is an extension
$$
1 \to \U \to \cG \to R \to 1
$$
where $\U$ is prounipotent.  Now we are in the situation of
Theorem~\ref{Th:unip-cohomology}. Denote the Lie algebra of $\U$ by $\u$. 
Since $\u$ is a $\cG$-module by the adjoint action, the natural weight
filtration on $\u$ induces one on $\Hcts^\dot(\u)$. By looking at cochains, it
is not difficult to see that if $\u = W_{-N}\u$ for some $N>0$, then 
\begin{equation}
\label{eq:Hmvanish}
W_n \Hcts^m(\u) = 0 \text{ if } n < Nm.
\end{equation}
Let $V_\alpha$ be  an irreducible $R$-module.  Since $w$ is central in $R$, the
$\Gm$-action commutes with the $R$-action, so Schur's Lemma implies that there
is an integer $n(\alpha)$ such that $\Gm$ acts on $V_\alpha$ via the
$n(\alpha)$th power of the standard character. This is the weight of $V_\alpha$
as a $\cG$-module. Now (\ref{eq:Hmvanish}) and Theorem~\ref{Th:unip-cohomology}
imply
$$
H^m(\cG,V_\alpha)=0
$$ 
if $n(\alpha) > -Nm$. Note that always $\u=W_{-1}\u$.

Suppose that $V$ is an $\l$-adic $\G$-module, i.e., a $\Ql$-vector space with
continuous $\G$-action. We shall need the continuous cohomology
$\Hcts^\dot(\G,V)$, which is defined as the cohomology of a suitable complex of
continuous cochains as in \cite[Sect.~2]{tate}.  A $\cG$-module $V$ can be
considered as an $\l$-adic $\G$-module through $\rhohat:\G \to \cG(\Ql)$. 
There is a natural group homomorphism
$$
\Phi^m: H^m(\cG,V) \to \Hcts^m(\G,V)
$$
for each $m\ge 0$.

Let $\{V_\alpha\}_\alpha$ be as in Theorem~\ref{Th:unip-cohomology}. These are
considered as $\G$-modules via $\rho$. The following theorem is our basic tool
for computing $\u$ when the appropriate continuous cohomology groups
$\Hcts^i(\G,V_\alpha)$ are known for $i=1,2$.

\begin{theorem}\label{h1h2}
For $m=1,2$, the mappings $\Phi^m$ defined above satisfy:
\begin{enumerate}
\item $\Phi^1:H^1(\cG,V_\alpha)\to \Hcts^1(\G,V_\alpha)$
is an isomorphism if $n(\alpha)<0$;
\item $\Phi^2:H^2(\cG,V_\alpha)\to \Hcts^2(\G,V_\alpha)$
is injective.
\end{enumerate}
\end{theorem}

This and Theorem~\ref{Th:unip-cohomology} imply the following, by using
(\ref{eq:Hmvanish}) and the comment following it.

\begin{corollary}
\label{cor:h1h2}
\begin{enumerate}
\item There is a natural $R$-equivariant isomorphism
$$
\Hcts^1(\u)\cong \bigoplus_{\{\alpha:n(\alpha)\leq -1\}}
\Hcts^1(\G,V_\alpha)\otimes V_\alpha^\ast.
$$
\item
If $N$ is an integer such that $\Hcts^1(\G,V_\alpha)=0$ for $0> n(\alpha)> -N$,
then there is a natural $R$-equivariant inclusion
$$
\Phi : \Hcts^2(\u) \hookrightarrow
\bigoplus_{\{\alpha:n(\alpha)\leq -2N\}}
\Hcts^2(\G,V_\alpha)\otimes V_\alpha^\ast.
$$
\end{enumerate}
\end{corollary}

This is proved in \cite{hain-matsumoto}. Below we shall give another more
categorical proof, similar to that in \cite{beilinson-deligne:orig}.

\begin{corollary}
\label{extras}
\begin{enumerate}
\item
If $\Hcts^1(\G,V_\alpha) = 0$ whenever $n(\alpha) < 0$, then $\u = 0$. \item
Let $N$ be as in Corollary~\ref{cor:h1h2}. If $\Hcts^2(\G,V_\alpha) = 0$
whenever $n(\alpha) \leq -2N$, then $\u$ is free as a pronilpotent Lie algebra.
\end{enumerate}
\qed
\end{corollary}

In the proof of Theorem~\ref{h1h2} we shall use Yoneda extensions. Let $V$ be a
finite dimensional $\l$-adic $\G$-module. For each $m \geq 0$, define
$$
\Ext^m_\G(\Ql,V)
$$
to be the $m$-th Yoneda extension group in the category of finite dimensional
$\l$-adic $\G$-modules, where $\Ql$ denotes the trivial $\G$-module. For each
$m\geq 1$, there is a natural homomorphism $$ \Ext^m_\G(\Ql, V) \to
\Hcts^m(\G,V), $$ which, by Theorem~\ref{th:cohomcomparison}, is an isomorphism
when $m=1$, and injective when $m=2$.

There is an exact functor from the category of weighted $\cG$-modules to
the category of $\l$-adic $\G$-modules.  It induces morphisms between the
extension groups, and hence homomorphisms
$$
\Psi^m: H^m(\cG,V) \to \Ext^m_\G(\Ql,V), \qquad m \ge 0.
$$
The homomorphisms $\Phi^m$ above factor through these:
$$
\xymatrix{
H^m(\cG,V) \ar[r]_{\Psi^m} \ar@/^1.5pc/[rr]^{\Phi^m} & \Ext^m_\G(\Ql,V) \ar[r] &
\Hcts^m(\G,V). \cr
}
$$
In fact, this is one of several equivalent ways to define the natural mappings
$\Phi^m$.

\begin{proof}[Proof of Theorem~\ref{h1h2}]
In view of Theorem~\ref{th:cohomcomparison}, it suffices to prove that $\Psi^1$
is an isomorphism, and that $\Psi^2$ is injective.

Since the functor from the category of weighted $\G$-modules to the category of
$\G$-modules is fully faithful,  a 1-step extension of weighted $\G$-modules
splits if it splits as an extension of $\G$-modules. This establishes the
injectivity of $\Psi^1$.

To prove surjectivity of $\Psi^1$, we define a natural weight filtration on each
$\G$-module extension $E$ of $\Ql$ by $V_\alpha$. Simply set $W_0 E =E$ and
$W_{-1}E=E$. Since $n(\alpha)<0$, this makes $E$ a weighted $\G$-module.

To prove that $\Psi^2$ is injective, we need to show that if a 2-step extension 
\begin{equation}
\label{eq:V2-step}
1 \to V_\alpha \to E_2 \to E_1 \to \Ql \to 1
\end{equation}
lies in the trivial class of extensions of $\G$-modules, then it also lies in
the trivial class of extensions of {\em weighted} $\G$-modules.

If $n(\alpha)\geq -1$, then $H^2(\cG,V_\alpha)=0$, and there is nothing to
prove. Thus we may assume $n(\alpha)\leq -2$. Since $W_m$ is an exact functor,
we may apply $W_0$ to (\ref{eq:V2-step}) to obtain  another 2-step extension,
without changing the extension class. Then, taking $\Gr^W_0$, we have a short
exact sequence
$$
0 \to \Gr^W_0 E_2 \to \Gr^W_0 E_1 \to \Ql \to 0
$$
of $R$-modules. Since $R$ is reductive, this has a splitting
$\Ql \hookrightarrow \Gr^W_0 E_1$. Taking the inverse images of this copy
of $\Ql$ along $E_2 \to E_1 \to \Gr^W_0 E_1$ in $E_2$ and in $E_1$, we obtain
a 2-step extension 
$$
0 \to V_\alpha \to E_2' \to E_1' \to \Ql \to 0
$$
equivialent to (\ref{eq:V2-step}) satisfying $W_{-1}E_2'=E_2'$ and
$W_{0}E_1'=E_1'$. Using the dual argument, we may assume that
(\ref{eq:V2-step}) satisfies $W_0E_1=E_1$, $W_{-1}E_2=E_2$,
$W_{n(\alpha)-1}E_2=0$, and $W_{n(\alpha)}E_1=0$.

By Yoneda's characterization \cite[p.~575]{yoneda2} of trivial $m$-step
extensions, the extension (\ref{eq:V2-step}) represents the trivial 2-step
extension class as $\G$-modules if and only if there is a $\G$-module $E$
and exact sequences
$$
0 \to E_2 \to E \to \Ql \to 0 \text{ and } 0 \to V_\alpha \to E \to E_1 \to 0
$$
which are compatible with the existing mappings $V_\alpha\hookrightarrow E_2$
and $E_1 \twoheadrightarrow \Ql$. To establish the injectivity of $\Psi^2$, it
suffices to prove that $E$ is a weighted $\G$-module. But $E$ has the weight
structure $W_0E=E$ and $W_{-1}E=E_2$. This completes the proof of
Theorem~\ref{h1h2}.
\end{proof}

\begin{example}
\label{zlx}
Suppose that $\G = \Zlx$, that $R=\Gm_{/\Ql}$ and that
$\rho : \Zlx \hookrightarrow \Gm(\Ql)=\Qlx$ is the natural inclusion.
Take $w$ to be the inverse of the square of the standard character. (With
this choice, representation theoretic weights coincide with the weights
from Hodge and Galois theory.) In this example we compute the
weighted completion of $\Zlx$ with respect to $\rho$ and $w$. Note that
$$
\Hcts^1(\Zlx,\Ql(n)) = 0,
$$
for all non-zero $n\in \Z$, where $\Ql(n)$ denotes the $n$th power of the
standard representation of $\Gm$. It has weight $-2n$ under the central
cocharacter.

Corollary~\ref{extras} tells us
that the unipotent radical $\U$ of the weighted completion of $\Zlx$ is
trivial, so that the weighted completion of $\Zlx$ with respect
to $\rho$ is just $\rho : \Zlx \to \Gm(\Ql)$. More generally, if $\G$ is an
open subgroup of $\Zlx$, then the weighted completion of $\G$, relative to
the restriction $\G \to \Qlx$ of the homomorphism $\rho$ above and the same
$w$, is simply $\Gm$.
\end{example}

\begin{example}
Let $\M_g$ be the moduli stack of genus $g$ curves over $\Spec\Z$. 
Suppose that there is a $\Z[1/\l]$-section $x:\Spec \Z[1/\l]\to \M_g$.
We allow tangential sections, and then such $x$ exist for all $g$.

Let $\bar{x}:\Spec \Qbar \to \Spec \Z[1/\l] \to \M_g$ be
a geometric point on the generic point of $x$. 
Let $C_{\bar{x}}$ be the curve corresponding to $\bar{x}$.

There is a short exact sequence of algebraic fundamental groups
\begin{equation}
\label{eq:alg-fund}
1 
\to \pi_1(\M_g\otimes\Qbar,\bar{x})
\to \pi_1(\M_g\otimes\Q,\bar{x})
\to G_\Q
\to 1,
\end{equation}
where the left group is isomorphic to the profinite completion $\Ghat_g$ of the
mapping class group $\G_g$ of a genus $g$ surface. We fix such an isomorphism.
We have the natural representation
\begin{equation}
\label{eq:moduli-rep}
\pi_1(\M_g\otimes\Q,\bar{x}) \to \Aut \Het^1(C_{\bar{x}},\Ql).
\end{equation}
It is known that the image of (\ref{eq:moduli-rep}) is isomorphic to
$\GSp_g(\Zl)$, where $\GSp_g$ denotes the group of symplectic similitudes of
a symplectic module of rank $2g$.

By considering the action of the mapping class group on the $\Z/\l\Z$ homology
of the surface, we obtain a natural representation $\Ghat_g \to \Sp_g(\Z/\l)$.
Let $\G_g^\l$ be the largest quotient of $\Ghat_g$ that also maps to
$\Sp_g(\Z/\l)$ and such that the kernel of the induced mapping
$\G_g^\l \to \Sp_g(\Z/\l)$ is a pro-$\l$ group.

One can construct a quotient
$$
1 \to \G_g^\l \to \G_g^{\arith,\l} \to G_{\Q,\{\l\}} \to 1
$$
of the short exact sequence (\ref{eq:alg-fund}) such that the homomorphism
(\ref{eq:moduli-rep}) induces a homomorphism
$$
\rho:\G_g^{\arith,\l} \to \GSp_g(\Ql)
$$
from (\ref{eq:moduli-rep}).

Define the central cocharacter $\w:\Gm \to \GSp_g$ by $x \to x^{-1}I_{2g}$. 
In \cite{hain-matsumoto:2} we show that the weighted completion
$\cG_g^{\arith,\l}$ of $\G_g^{\arith,\l}$ is an extension
$$
\cG_g\otimes \Ql \to \cG_g^{\arith,\l} \to \A^\l_{\Q,\{\l\}} \to 1,
$$
where $\cG_g$ is the completion of $\G_g$ relative to the standard homomorphism
$\rho : \G_g \to \Sp_g(\Q)$, which is studied in \cite{hain:comp} and for which
a presentation is given in \cite{hain:torelli}. We expect
that the left homomorphism is injective.
\end{example}

\section{Computation of $\AFSl$}
\label{compn}

In this section, we compute $\AFSl$, the tannakian fundamental group of the
category of $\l$-adic mixed Tate modules over $\XFS$. An equivalent computation
was done by Beilinson and Deligne in \cite{beilinson-deligne:orig}. We shall
need the following result of Soul\'e \cite{soule} when $\l$ is odd. The case
$\l=2$ follows from \cite{Rognes}. Recall that $d_n$ is defined in (\ref{dn}).

\begin{theorem}[Soul\'e \cite{soule}]
\label{Th:Galcoh}
With notation as above,
$$
K_{2n-1}(\XFS)\otimes \Ql \cong \Hcts^1(\GFS,\Ql(n))
$$
and hence
$$
\dim_\Ql \Hcts^1(\GFS,\Ql(n)) =
\begin{cases}
d_1 +\#S & n = 1, \cr
d_n & n > 1.
\end{cases}
$$
In addition, $\Hcts^2(\GFS,\Ql(n))$ vanishes for all $n\ge 2$.
\end{theorem}

Denote the unipotent radical of $\AFSl$ by $\KFSl$. We have the exact
sequence
$$
1 \to \KFSl \to \AFSl \to \Gm \to 1,
$$
and the corresponding exact sequence of Lie algebras
$$
0 \to \kFSl \to \aFSl \to \Ql \to 0.
$$

The Lie algebra $\aFSl$, being the Lie algebra of a weighted completion,
has a natural weight filtration. Note that since $w$ is the inverse 
of the square of the standard character, all weights are even. Thus the
weight filtration of $\aFSl$ satisfies
$$
\aFSl = W_0 \aFSl, \quad \kFSl= W_{-2}\aFSl
\text{ and } \Gr^W_{2n+1} \aFSl = 0 \text{ for all }n\in \Z,
$$
and we may take $N=2$ in Corollary~\ref{cor:h1h2}.
The basic structure of $\AFSl$ now follows from Corollary~\ref{cor:h1h2},
Corollary~\ref{extras} and Soul\'e's computation above.

\begin{theorem}[Hain-Matsumoto \cite{hain-matsumoto}]
\label{main}
\label{a_l}
The Lie algebra $\Gr^W_\dot\kFSl$ is a free Lie algebra and there is a natural
$\Gm$-equivariant isomorphism
$$
\Hcts^1(\kFSl) \cong \bigoplus_{n=1}^\infty\Hcts^1(\GFS,\Ql(n))\otimes\Ql(-n)
\cong \Ql(-1)^{d_1 + \#S} \oplus \bigoplus_{n>1} \Ql(-n)^{d_n},
$$
where $d_n$ is defined in (\ref{dn}). Any lift of a graded basis of
$H_1(\Gr^W_\dot \kFSl)$ to a graded set of elements of $\Gr^W_\dot \kFSl$
freely generates $\Gr^W_\dot \kFSl$.
\end{theorem}

As a corollary of the proof, we have:

\begin{corollary}
There are natural isomorphisms
$$
\Ext^m_{\MTMlFS}(\Ql,\Ql(n))
\cong
\begin{cases}
\Ql & \text{ when $m = n = 0$,} \cr
\Hcts^1(\GFS,\Ql(n)) & \text{ when $m=1$ and $n>0$,} \cr
0  & \text{ otherwise.}
\end{cases}
$$
Consequently, for all $n\in \Z$, there are natural isomorphisms
$$
\Ext^1_{\MTMlFS}(\Ql,\Ql(n)) \cong K_{2n-1}(\Spec\OFS)\otimes\Ql.
$$
\end{corollary}

\begin{corollary}\label{cor:post2-del}
Suppose that there is a category of mixed Tate motives $\T(\XFS)$ with
properties (i)--(v) as in Section~\ref{sec:mtm}. If Deligne's
conjecture~\ref{del_conj} is true, then, the image of the  $\l$-adic
realization functor $\real_\l$ in (\ref{realization}) is equivalent to the
category of weighted $\GFS$-modules. In particular, Postulate~\ref{post2}
follows. 
\end{corollary}

\begin{proof}
Deligne's conjecture~\ref{del_conj} implies that $\pi_1(\T(\XFS))$, the
tannakian fundamental group, is an extension
$$
1 \to \U_\XFS \to \pi_1(\T(\XFS)) \to \Gm \to 1,
$$
where $\U_\XFS$ is a free prounipotent group  generated by
$K_{2n-1}(\XFS)^\ast$. This and Theorem~\ref{main} show that the natural map
$\pi_1(\T_\l(\XFS)) \to \pi_1(\T(\XFS))\otimes \Ql$ is an isomorphism, and it
follows that
$$
\real_\l : \T(\XFS)\otimes \Ql \to \text{$\l$-adic $G_F$-modules}
$$
is fully faithful and its image is equivalent to the category of weighted
$\GFS$-modules.
\end{proof}

Note that these theorems can be generalized to the 
case where $S$ may not contain all the primes above $\l$,
see Section~\ref{sec:crys}.

\subsection*{Another Example}
Suppose that $S$ is a finite set of  rational primes containing $\l$. Suppose
that $F$ is a  finite Galois extension of $\Q$ with Galois group $G$, which is
unramified outside $S$. Define
$$
\rho : G_{\Q,S} \to \Gm(\Ql) \times G
$$
by
$$
\rho(\sigma) = (\chi_\l(\sigma), f(\sigma))
$$
where $f : G_{\Q,S} \to G$ is the quotient homomorphism and $\chi_\l$
is the $\l$-adic cyclotomic character. Define
$$
w : \Gm \to \Gm \times G
$$
by $w : x \mapsto (x^{-2},1)$. It is a central cocharacter. Denote the
weighted completion of $G_{\Q,S}$ with respect to $\rho$ and $w$ by
$\cG_{\Q,S}$.

Denote the set of primes in $\O_F$ that lie over $S \subset \Spec\Z$ by $T$.

\begin{proposition}
There is a natural inclusion $\iota : \A_{F,T}^\l \to \cG_{\Q,S}$ and an
exact sequence
$$
\begin{CD}
1 @>>> \A_{F,T}^\l @>\iota>> \cG_{\Q,S} @>>> G @>>> 1.
\end{CD}
$$
\end{proposition}

\begin{proof}
If $\{V_\alpha\}$ is a set of representatives of the isomorphism classes of
irreducible representations of $G$, then $\{\Ql(m)\boxtimes V_\alpha\}$ is a
set of representatives of the isomorphism classes of irreducible
representations of $\Gm \times G$, where $W \boxtimes V$ denotes the exterior
tensor product of a representation $W$ of $\Gm$ and $V$ of $G$. 
Consider the restriction mapping
$$
\phi: \Hcts^i(G_{\Q,S},\Ql(m)\boxtimes V_\alpha)
\to
\Hcts^i(G_{F,T},\Ql(m)\boxtimes V_\alpha)^G
$$
and the transfer mapping \cite{tate}
$$
\psi: \Hcts^i(G_{F,T},\Ql(m)\boxtimes V_\alpha)^G
\to
\Hcts^i(G_{\Q,S},\Ql(m)\boxtimes V_\alpha).
$$
A direct computation on cocycles shows that $\phi \circ \psi$ and $\psi \circ
\phi$ are both multiplication by the order of $G$, and are thus isomorphisms.

So $\Hcts^i(G_{\Q,S},\Ql(m)\boxtimes V_\alpha)$ vanishes if $V_\alpha$ is
non-trivial, and is $\Hcts^i(G_{F,T},\Ql(m))$ if $V_\alpha$ is trivial.  This
shows that the unipotent radical of  the completion  $\cG_{\Q,S}$ is isomorphic
to that of $\A_{F,T}^\l$.

By functoriality of weighted completion, we have a homomorphism $\A_{F,T}^\l
\to \cG_{\Q,S}$ which induces the isomorphism on the unipotent radical. The
statement follows.
\end{proof}

\section{Applications to Galois Actions on Fundamental Groups}
\label{proof}

Let $\Gl$ denote $G_{\Q,\{\l\}}$. In this section, we sketch how our
computation of the weighted completion of $\Gl$ can be used to prove Deligne's
Conjecture~\ref{del-ihara_conj} about the action of the absolute Galois group
$G_\Q$ on the pro-$\l$ completion of the fundamental group of $\PminusC$.
Modulo a few technical details, which are addressed in \cite{hain-matsumoto},
the proof proceeds along the expected lines suggested in
Section~\ref{discussion} given the computation of $\AFSl$.

We begin in a more general setting. Suppose that $F$ is a number field and that
$X$ is a variety over $F$. Set $\Xbar = X \otimes \Qbar$ and denote the
absolute Galois group of $F$ by $G_F$. Suppose that the \'etale cohomology
group $\Het^1(\Xbar,\Ql(1))$ is a trivial $G_F$-module. Let $S$ be a set of
finite primes of $F$, containing those above $\l$. Suppose that $X$ has a model
$\X$ over $\Spec\OFS$ which has a base point section $x: \Spec\OFS \to \X$
(possibly tangential) such that $(\X,x)$ has good reduction outside
$S$.\footnote{What we mean here is that there is a scheme $\cXtilde$, proper
over $\Spec\OFS$, and a divisor $D$ in $\cXtilde$ which is relatively normal
crossing over  $\Spec\OFS$ such that  $\X = \cXtilde - D$, and $D$ does not
intersect with $x$.  In the tangential case,  the tangent vector should be
non-zero over each point of $\Spec\OFS$.} Then the $G_F$-action on the pro-$\l$
fundamental group $\pi_1^\l(\Xbar,x)$ factors through $\GFS$.

Denote the $\l$-adic unipotent completion of $\pi_1^\l(\Xbar,x)$ by $\cP$
(see \cite[Appendix~A]{hain-matsumoto}) and its Lie algebra by $\p$.
The lower central series filtration of $\p$ gives it the structure of a
pro-object of the category $\T_\l(\XFS)$ of $\l$-adic mixed Tate modules
over $\XFS$. It follows that the $G_F$-action on $\cP$ induces a homomorphism
$$
\AFSl \to \Aut\cP \cong \Aut \p
$$
and that the action of $G_F$ on $\cP$ factors through the composition of
this with the natural homomorphism $G_F \to \GFS \to \AFSl(\Ql)$.
One can show (see \cite[Sect.~8]{hain-matsumoto}) that the image of
$G_{F(\mu_{\l^\infty})}$ in $\AFSl$ lies in and is Zariski dense in $\KFSl$.

For the remainder of this section, we consider the case where $X = \Pminus$,
$F=\Q$, $S=\{\l\}$ and $x$ is the tangential base point $\v$. Goncharov's
conjecture \cite[Conj.~2.1]{goncharov} (cf.\ the generation part of
Conjecture~\ref{gonch_conj}) follows immediately, since $\k_{\Q,\{\l\}}$ is
generated by $z_1,z_3,z_5,\ldots$, where $z_j$ has weight $-2j$. The image of
$z_1$ can be shown to be trivial.

We are now ready to give a brief sketch of the proof of
Conjecture~\ref{del-ihara_conj}. One can define a filtration $\cI_\l^\dot$ on
$G_\Q$ similar to $I_\l^\dot$ using the lower central series $L^\dot\cP$ of
$\cP$ instead:
$$
\cI_\l^m G_\Q = \ker\{G_\Q \to \Out \cP/L^{m+1}\cP\}
$$
where $L^m\cP$ is the $m$th term of its lower central series. The lower
central series of $\cP$ is related to its weight filtration by
$$
W_{-2m} \cP = L^m\cP, \quad \Gr^W_{2m+1}\cP = 0.
$$
There is a natural isomorphism (see \cite[Sect.~10]{hain-matsumoto})
$$
[\Gr_\I^m G_\Q]\otimes \Ql \cong [\Gr_\cIell^m G_\Q]\otimes \Ql.
$$
Thus it suffices to prove that $[\Gr_\cIell^m G_\Q]\otimes \Ql$ is generated
by elements $s_3, s_5, s_7, \dots $, where $s_j$ has weight $-2j$.

As above, the homomorphism $G_\Q \to \Out\cP$ factors through the sequence
$$
G_\Q \to G_{\Q,\{\l\}} \to \A_{\Q,\{\l\}} \to \Out \cP
$$
of natural homomorphisms. A key point (\cite[Sect.~8]{hain-matsumoto}) is that
the image of $\I^1G_\Q$ in $\K_{\Q,\l}^\l$ is Zariski dense.  This and the
strictness can be used to establish isomorphisms
\begin{multline*}
Gr_\I^m G_\Q \otimes \Ql \cong \Gr^W_{-2m} (\im\{\k_{\Q,\l}^\l\to \OutDer\p\})
\cr
\cong \im\{\Gr^W_{-2m}\k_{\Q,\l}^\l \to \Gr^W_{-2m}\OutDer\p\}
\end{multline*}
for each $m>0$.

Theorem~\ref{main} implies that $\Gr^W_\dot \k_{\Q,\l}^\l$ is freely
generated by $\sigma_1, \sigma_3, \sigma_5, \ldots$
where $\sigma_{2i+1} \in \Gr^W_{-2(2i+1)}\k_{\Q,\l}^\l$. It is easy to show
that the image of $\sigma_1$ vanishes in $\Gr^W_\dot\OutDer\p$. It follows
that the image of $\Gr^W_\dot\k_{\Q,\l}^\l$ is generated by the images of
$\sigma_3, \sigma_5, \sigma_7, \ldots$, which completes the proof.

\begin{remark}
Ihara proves the openness of the group generated by $\sigma_{2i+1}$ in a
suitable Galois group, see \cite{Ihara9}. He also establishes the non-vanishing
of the images of the $\sigma_{2i+1}$ and some of their brackets in
\cite{Ihara3}.
\end{remark}

\section{When $\l$ is not contained in $S$}\label{sec:crys}

Let $[\l]$ denote the set of all primes above $\l$ in $\O_F$.  In this section,
we generalize the definition of the category $\T_\l(\XFS)$ of $\l$-adic mixed
Tate modules smooth over $\XFS=\Spec\O_F-S$ (see Section~\ref{wtd_ell_mods}) to
the case where $S$ does not necessarily contain $[\l]$.

For this,  we define the category $\T_\l(\XFS)$ of $\l$-adic mixed Tate modules
over $\XFS$ to be the full subcategory  of $\T_\l(X_{F, S\cup [\l]})$  (defined
in Section~\ref{wtd_ell_mods}) consisting of the Galois modules  which are
crystalline at every prime $\p \in [\l] - S$. (Recall that an $\l$-adic
$G_F$-module $M$ is crystalline at a prime $\p$ of $F$ if it is crystalline as 
$G_{F_\p}$-module, where $F_\p$ is the completion of $F$ at $\p$ and $G_{F_\p}$
is identified with the decomposition group of $G_F$ at $\p$, see \cite{FO,BK}
for crystalline representations.)

It is known that the crystalline property is closed under tensor products,
direct sums, duals, and subquotients \cite{FO}, so that $\T_\l(\XFS)$ is a
tannakian category. Denote its tannakian fundamental group by $\AFSl$. We have
a short exact sequence
$$
1 \to \KFSl \to \AFSl \to \Gm \to 1,
$$
and the corresponding exact sequence of Lie algebras
$$
0 \to \kFSl \to \aFSl \to \Ql \to 0.
$$

Let $V$ be a $\GFS$-module. The {\it finite part} of the first degree Galois
cohomology $\Hctsf^1(\GFS, V)\subset \Hcts^1(\GFS,V)$  is defined in
\cite[(3.7.2)]{BK}. This corresponds to those extensions of $\Ql$ by $V$ as
$\GFS$-modules, which are crystalline at every prime in $[\l]$ outside $S$. By
a remark on p.~354 in \cite{BK}, $\Hctsf^1(\GFS,\Ql)=(\OFS^\times)\otimes_\Zl
\Ql$, so its dimension is $d_1 + \#S = r_1 + r_2 + \#S - 1$. 
Theorem~\ref{main} is generalized as follows, by replacing $\Hcts^1$ with
$\Hctsf^1$ \cite{hain-matsumoto}. We shall give a categorical proof below.

\begin{theorem}
\label{maincrys}
The Lie algebra $\Gr^W_\dot\kFSl$ is a free Lie algebra and there is a natural
$\Gm$-equivariant isomorphism
$$
\Hcts^1(\kFSl) \cong \bigoplus_{n=1}^\infty\Hctsf^1(\GFS,\Ql(n))\otimes\Ql(-n)
\cong \Ql(-1)^{d_1 + \#S} \oplus \bigoplus_{n>1} \Ql(-n)^{d_n},
$$
where $d_n$ is defined in (\ref{dn}). Any lift of a graded basis of
$H_1(\Gr^W_\dot \kFSl)$ to a graded set of elements of $\Gr^W_\dot \kFSl$
freely generates $\Gr^W_\dot \kFSl$.
\end{theorem}

\begin{corollary}
There are natural isomorphisms
$$
\Ext^m_{\MTMlFS}(\Ql,\Ql(n))
\cong
\begin{cases}
\Ql & \text{ when $m = n = 0$,} \cr
\Hctsf^1(\GFS,\Ql(n)) & \text{ when $m=1$ and $n>0$,} \cr
0  & \text{ otherwise.}
\end{cases}
$$
Consequently, for all $n\in \Z$, there are natural isomorphisms
$$
\Ext^1_{\MTMlFS}(\Ql,\Ql(n)) \cong K_{2n-1}(\Spec\OFS)\otimes\Ql.
\qed
$$
\end{corollary}

This shows that $\T_\l(\XFS)$ has all the properties of the category
$\T(\XFS)\otimes \Ql$, where $\T(\XFS)$ is the category whose existence is
conjectured by Deligne. In particular,  $\Gr^W_\dot \k_{\Q,\emptyset}^\l$ is
free with generators $\sigma_3,\sigma_5,\ldots$.

\begin{proof}[Proof of Theorem~\ref{maincrys}]
It suffices to show that the natural mapping
$$
\Phi^1:H^1(\AFSl,\Ql(n))\to \Hctsf^1(\GFS,\Ql(n))
$$
is an isomorphism when $n\geq 1$ and that the natural mapping
$$
\Phi^2:H^2(\AFSl,\Ql(n))\to {\Hcts^2}(\GFS,\Ql(n))
$$
is injective when $n\ge 2$S.  The proof is similar to that of
Theorem~\ref{h1h2}. To show that $\Phi^1$ is an isomorphism, it suffices to
show that an extension $E$ of $\Ql$ by $\Ql(n)$ corresponding to an element of
$\Hctsf^1(\GFS,\Ql(n))$ is crystalline, which is well-known. So the first
assertion follows.

We now consider the case of $\Phi^2$. Set $V_\alpha = \Ql(n)$. We may assume
$n\geq 2$.  It suffices to show that $E$ in the proof of Theorem~\ref{h1h2} is
crystalline provided $E_1$ and $E_2$  are crystalline. But this follows from
the next result, which will be proved below.

\begin{proposition}\label{prop:crys-ext}
Let 
$$
0\longrightarrow V\longrightarrow U \longrightarrow 
\Ql(1)^{n}\longrightarrow 0
$$
be a short exact sequence of crystalline $\l$-adic representations of
$G_{F_\p}$.  Assume that $V$ is a successive extension of direct sums of 
a finite number of copies of $\Ql(r)$ with $r\geq 2$. Then, for any extension
$$
0\longrightarrow U\longrightarrow E\longrightarrow \Ql\longrightarrow 0
$$
of $\l$-adic representations of $G_{F_\p}$,  $E$ is crystalline if and only if
its pushout by the surjection $U\to \Ql(1)^n$ is crystalline.
\end{proposition}

Let $U$ be $E_2$ as in the proof of Theorem~\ref{Th:unip-cohomology}. Since
$W_{-2}E_2=E_2$, $U$ is an extension of $\Ql(1)^n$ for some $n$. Since $m\geq
2$,  the pushout of $E$ along $U \to \Ql(1)^n$ is a quotient of $E_1$, and
hence is crystalline. Thus the proposition says that $E$ is crystalline, which
completes the proof of Theorem~\ref{maincrys}.
\end{proof}

Proposition~\ref{prop:crys-ext} follows from the 
following two lemmas.
\begin{lemma}\label{lem:1}
Let 
$$
0\to V_1 \to V_2 \to V_3 \to 0
$$
be a short exact sequence of crystalline $\l$-adic
representations of $G_{F_\p}$. Then we have a long exact sequence
\begin{eqnarray*}
0&\to& H^0(G_{F_\p},V_1)\to H^0(G_{F_\p},V_2)
\to H^0(G_{F_\p},V_3) \\ 
&&
\to \Hctsf^1(G_{F_\p},V_1)
\to \Hctsf^1(G_{F_\p},V_2)
\to \Hctsf^1(G_{F_\p},V_3)
\to 0. \qed
\end{eqnarray*}
\end{lemma}
This follows from \cite[Cor. 3.8.4]{BK}.

\begin{lemma}\label{lem:2}
Let $V$ be a crystalline $\l$-adic representation of 
$G_{F_\p}$. 
If $V$ is a successive extension of $\Ql(r)$ $(r\geq 2)$, then
$\Hctsf^1(G_{F_\p},V)=\Hcts^1(G_{F_\p},V)$.
\end{lemma}

\begin{proof}
The proof is by induction on the dimension of $V$. In the  case dim$(V)=1$,
this is well-known (loc. cit. Example 3.9). Assume $\dim V=n \geq 2$ and the
claim is true for $n-1$. By assumption, there exists an exact sequence of
$\l$-adic representations of $G_{F_\p}$:
$$
0\longrightarrow V'\longrightarrow V\longrightarrow 
\Ql(r)\longrightarrow 0
$$
for some integer $r\geq 2$ such that $V'$ satisfies the assumption 
of the lemma. By Lemma~\ref{lem:1}, we have the following commutative diagram
whose two rows are exact:
{\scriptsize{
$$
\begin{CD}
0 @>>> \Hctsf^1(G_{F_\p},V')@>>> 
\Hctsf^1(G_{F_\p},V) @>>> \Hctsf^1(G_{F_\p},\Ql(r)) @>>>0 \\
@. @VVV  @VVV @VVV @.\\
0 @>>> \Hcts^1(G_{F_\p},V') @>>> \Hcts^1(G_{F_\p},V) @>>> 
\Hcts^1(G_{F_\p},\Ql(r))
\end{CD}
$$
}}
The right vertical arrow is an isomorphism and the left one is also
an isomorphism by the induction hypothesis. Hence the middle one is
also an isomorphism. 
\end{proof}

\begin{proof}[Proof of Proposition~\ref{prop:crys-ext}]
By Lemma~\ref{lem:1}, we  have the following commutative diagram whose two rows
are exact:
{\scriptsize{
$$
\CD
0 @>>> \Hctsf^1(G_{F_\p},V) @>>> \Hctsf^1(G_{F_\p},U) @>>> 
\Hctsf^1(G_{F_\p},\Ql(1)^{n}) @>>> 0\\
@. @VVV @VVV @VVV\\
0 @>>> \Hcts^1(G_{F_\p},V) @>>> \Hcts^1(G_{F_\p},U) @>>> 
\Hcts^1(G_{F_\p},\Ql(1)^{n})
\endCD
$$
}}
and the left vertical arrow is an isomorphism by Lemma~\ref{lem:2}. Hence
the right square is cartesian.
\end{proof}

\appendix

\section{Continuous Cohomology and Yoneda Extensions}

In this appendix we prove a result about the relation between continuous
cohomology and Yoneda extension groups in low degrees. It is surely well
known, but we know of no reference.

Suppose that $K$ is a topological field, and $\G$ a topological group. A {\em
continuous $\G$-module} is a $\G$-module $V$, where $V$ is a finite dimensional
$K$-vector space. The action $\G \to \GL(V)$ is required to be continuous,
where $\GL(V)$ is given the topology induced from that of $K$.

Denote by $\cC(\G,K)$ the category of  finite dimensional continuous
$\G$-modules.  Since any $K$-linear morphism between finite dimensional vector
spaces is continuous, this is an abelian category. For continuous $\G$-modules
$A$ and $B$, define $\Ext_{\G}^\dot(A,B)$ to be the graded group of Yoneda
extensions of $B$ by $A$ in the category $\cC(\G,K)$.

For a continuous $\G$-module $A$, one also has the continuous cohomology groups
$\Hcts^\dot(\G,A)$ defined by Tate \cite{tate}, which are defined using the
complex of continuous cochains.

\begin{theorem}\label{th:cohomcomparison}
If $A$ is a continuous $\G$-module, then there is a natural isomorphism
$\Ext_{\G}^1(K,A) \cong \Hcts^1(\G,A)$ and a natural injection
$\Ext_{\G}^2(K,A) \hookrightarrow \Hcts^2(\G,A)$.
\end{theorem}

\begin{proof}
It is well known that an extension $0 \to A \to E \to K \to 0$ in $\cC(\G,K)$
gives a continuous cocycle $f : \G \to A$ by choosing a lift $e \in E$ of
$1\in K$ and defining $f(\sigma) = \sigma(e)-e$. Conversely, for a given
continuous cocycle $f$, we may define continuous $\G$-action on $A \oplus K$ by
$\sigma : (a,k) \mapsto (\sigma(a)+kf(\sigma), k)$. These are mutually inverse,
which establishes the first claim.

To prove the second claim, we first define a $K$-linear mapping
$$
\varphi : \Ext_\G^2(K,A) \to \Hcts^2(\G,A)
$$
as follows. For $c\in \Ext_\G^2(K,A)$, choose a 2-fold extension
$0 \to A \to E_2 \to E_1 \to K \to 0$ that represents it. By \cite{yoneda2},
$c$ is the image under the connecting homomorphism
$$
\delta : \Ext_\G^1(K,E_2/A) \to \Ext_\G^2(K,A)
$$
of the class $\tilde{c}$ of the extension $0 \to E_2/A \to E_1 \to K \to 0$.

We shall construct $\varphi$ so that the diagram
$$
\xymatrix{
\Ext_\G^1(K, E_2) \ar[r] \ar[d]^\simeq &
\Ext_\G^1(K, E_2/A) \ar[r]^\delta \ar[d]^\simeq_\psi &
\Ext_\G^2(K, A) \ar@{.>}[d]^\varphi \cr
\Hcts^1(\G, E_2)\ar[r] & \Hcts^1(\G, E_2/A) \ar[r]^{\delta_\cts}
& \Hcts^2(\G, A) \cr
}
$$
commutes, where the rows are parts of the standard long exact sequences
constructed in \cite{yoneda2} and \cite[Sect.~2]{tate}. Define $\varphi(c)$ to
be $\delta_\cts(\psi(\tilde{c}))$.

To prove $\varphi(c)$ is well-defined, it suffices to show that two 2-fold
extensions that fit into a commutative diagram
$$
\begin{CD}
0 @>>> A @>>> E_2' @>>>  E_1' @>>> K @>>> 0 \cr
@. @| @VVV @VVV  @| \cr
0 @>>> A @>>> E_2  @>>>  E_1  @>>> K @>>> 0 \cr
\end{CD}
$$ 
give a same element of $\Hcts^2(\G,A)$. But this follows from the functoriality of
the connecting homomorphism for $\Hcts^\dot$, i.e., the commutativity of
$$
\begin{CD}
\Hcts^1(\G,E_2'/A) @>>> \Hcts^2(\G,A) \cr
@VVV  @| \cr
\Hcts^1(\G,E_2 /A) @>>>  \Hcts^2(\G,A). \cr
\end{CD}
$$

The $K$-linearity of $\varphi$ is easily checked. Finally, the injectivity of
$\varphi$ follows from the fact that for each extension as above, $\varphi$ is
injective on the image of the connecting homomorphism
$\delta : \Ext_\G^1(K, E_2/A) \to \Ext_\G^2(K, A)$.
\end{proof}

Note that one may define
$$
\Ext_\G^m(K,A) \to \Hcts^m(G,A)
$$
by induction on $m$ in the same way.

\medskip

\noindent{\it Acknowledgements:} We would like to thank Marc Levine for
clarifying several points about motivic cohomology and Owen Patashnick for his
helpful comments on the manuscript.  We are indebted to Kazuya Kato and Akio
Tamagawa for pointing out a subtlety regarding continuous cohomology related
to  Theorem~\ref{th:cohomcomparison},  and to Takeshi Tsuji for the proof of
Proposition~\ref{prop:crys-ext}. We would also like to thank Sasha Goncharov
for pointing out the existence and relevance of the unpublished  manuscript
\cite{beilinson-deligne:orig} of Beilinson and Deligne. Finally, we would like
to thank the referee for doing a very thorough job and for many useful
comments.


\begin{thebibliography}{99}

\bibitem{beilinson}
A.~Beilinson:
{\it Higher regulators and values of $L$-functions} (Russian), Current problems
in mathematics, Vol. 24 (1984), 181--238.

\bibitem{beilinson-deligne:orig}
{A.~Beilinson, P.~Deligne:}
{\it Motivic Polylogarithms and Zagier's Conjecture}, unpublished
manuscript, 1992.

\bibitem{beilinson-deligne}
A.~Beilinson, P.~Deligne:
{\it Interpr\'etation motivique de la conjecture de Zagier reliant
polylogarithmes et r\'egulateurs}, Motives (Seattle, WA, 1991), Proc.\
Sympos.\ Pure Math., 55, Part 2, Amer.\ Math.\ Soc., Providence, RI, 1994,
97--121.

\bibitem{bms}
A.~Beilinson, R.~MacPherson, V.~Schechtman:
{\it Notes on motivic cohomology}, Duke Math.\ J.\ 54 (1987), 679--710.
 
\bibitem{bloch}
S.~Bloch:
{\it Algebraic cycles and higher $K$-theory}, Adv.\ in Math.\ 61 (1986),
267--304.

\bibitem{bloch:corr}
S.~Bloch:
{\it The moving lemma for higher Chow groups}, J.\ Algebraic Geom.\ 3 (1994),
537--568.

\bibitem{BK}
S.~Bloch, K.~Kato:
{\it L-Functions and Tamagawa Numbers of Motives},
The Grothendieck Festschrift Volume I,
Progress in Math. Vol.86, Birkh\"auser, (1990), 333-400.

\bibitem{bloch-kriz}
S.~Bloch, I.~Kriz:
{\it Mixed Tate motives}, Ann.\ of Math.\ 140 (1994), 557--605.

\bibitem{bloch-ogus}
S.~Bloch, A.~Ogus:
{\it Gersten's conjecture and the homology of schemes}, Ann.\ Sci.\ \'Ecole
Norm.\ Sup.\ (4) (1974), 181--201.

\bibitem{borel}
A.~Borel:
{\it Stable real cohomology of arithmetic groups}, Ann.\ Sci.\ \'Ecole Norm.\
Sup.\ (4) 7 (1974), 235--272 (1975).

\bibitem{borel:reg}
A.~Borel:
{\it Cohomologie de ${\rm SL}_{n}$ et valeurs de fonctions z\^{e}ta aux points
entiers}, Ann.\ Scuola Norm.\ Sup.\ Pisa Cl.\ Sci.\ (4) 4 (1977), 613--636.

\bibitem{borel-serre}
A.~Borel, J.-P.~Serre:
{\it Le th\'eor\`eme de Riemann-Roch}, Bull.\ Soc.\ Math.\ France 86
1958, 97--136.

\bibitem{dmos}
P.~Deligne and  J.~Milne:
{\it Tannakian categories}, in {Hodge Cycles, Motives, and Shimura Varieties}, 
(P.~Deligne, J.~Milne, A.~Ogus, K.-Y.~Shih editors), Lecture Notes in
Mathematics 900, Springer-Verlag, 1982.

\bibitem{deligne}
P.~Deligne:
{\it Le groupe fondamental de la droite projective moins trois points},
in Galois groups over $\Q$ (Berkeley, CA, 1987), 79--297, Math.\ Sci.\
Res.\ Inst.\ Publ., 16, Springer, New York-Berlin, 1989.

\bibitem{drinfeld}
V.~Drinfeld:
{\it On quasitriangular quasi-Hopf algebras and on a group that is
closely connected with $\Gal\Qbar/\Q)$}, (Russian) Algebra i Analiz 2 (1990),
149--181; translation in Leningrad Math.\ J.\ 2 (1991), no. 4, 829--860.

\bibitem{FO}
J.M.~Fontaine:
{\it Sur certains types de repr\'esentations $p$-adiques du groupe de Galois
d'un corps local: construction d'un anneau de Barsotti-Tate}, Ann.\ of Math.\
115 (1982), 529--577.


\bibitem{gillet}
H.~Gillet:
{\it Riemann-Roch theorems for higher algebraic $K$-theory}, Adv.\ in Math.\
40 (1981), 203--289.

\bibitem{goncharov:mtm}
A.~Goncharov:
{\it Volumes of hyperbolic manifolds and mixed Tate motives}, J.\ Amer.\
Math.\ Soc.\ 12 (1999), 569--618.

\bibitem{goncharov}
A.~Goncharov:
{\it Multiple zeta values, Galois groups, and geometry of modular
varieties}, preprint 2000.

\bibitem{hain:polylogs}
R.~Hain:
{\it Classical polylogarithms}, Motives (Seattle, WA, 1991), Proc.\ Sympos.\
Pure Math., 55, Part 2, Amer. Math. Soc., Providence, RI, 1994, 3--42.

\bibitem{hain:comp}
R.~Hain:
{\it Completions of mapping class groups and the cycle $C-C^-$},
in {\it Mapping Class Groups and Moduli Spaces of Riemann Surfaces},
Contemp.\ Math.\ 150 (1993), 75--105.

\bibitem{hain:torelli}
R.~Hain:
{\it Infinitesimal presentations of the Torelli groups}, J.\ Amer.\ Math.\
Soc.\ 10 (1997), 597--651.

\bibitem{hain-matsumoto}
R.~Hain, M.~Matsumoto:
{\it Weighted Completion of Galois Groups and Galois Actions on Fundamental
Groups}, Compositio Math., to appear, {\sf math.AG/0006158}.

\bibitem{hain-matsumoto:2}
R.~Hain, M.~Matsumoto:
{\it Completions of Arithmetic Mapping Class Groups}, in preparation.

\bibitem{hiller}
H.~Hiller:
{\it $\lambda $-rings and algebraic $K$-theory}, J.\ Pure Appl.\ Algebra 20
(1981), 241--266.

\bibitem{Ihara1}
Y.~Ihara:
{\it Profinite braid groups, Galois representations
and complex multiplications}, Ann.\ of Math., 123 (1986), 43--106.

\bibitem{Ihara3}
Y.~Ihara:
{\it The Galois representation arising from 
$\P^1-\{0,1,\infty\}$ and Tate twists of even degree},
in {\it Galois groups over $\Q$}, Publ.\ MSRI,
No.~16 (1989), Springer-Verlag, 299--313.

\bibitem{Ihara9}
Y.~Ihara:
{\it Some arithmetic aspects of Galois actions on
the pro-$p$ fundamental group of $\P^1-\{0,1,\infty\}$},
RIMS preprint 1229, 1999.

\bibitem{jannsen}
U.~Jannsen:
{\it Mixed motives and algebraic $K$-theory},
Lecture Notes in Mathematics, 1400, Springer-Verlag, Berlin, 1990.

\bibitem{Jantzen}
J.~Jantzen:
{\it Representations of Algebraic Groups},
Pure and Applied Mathematics Vol.131, Academic Press, 1987.

\bibitem{kratzer}
C.~Kratzer:
{\it $\lambda $-structure en $K$-th\'eorie alg\'ebrique}, Comment.\ Math.\
Helv.\ 55 (1980), 233--254.

\bibitem{levine:tate}
M.~Levine:
{\it Tate motives and the vanishing conjectures for algebraic $K$-theory},
Algebraic $K$-theory and algebraic topology (Lake Louise, AB, 1991), NATO
Adv.\ Sci.\ Inst.\ Ser.\ C Math.\ Phys.\ Sci., 407, Kluwer Acad. Publ.,
Dordrecht, 1993, 167--188.

\bibitem{levine}
M.~Levine:
{\it Bloch's higher Chow groups revisited}, $K$-theory (Strasbourg, 1992),
Ast\'erisque No.\ 226, (1994), 10, 235--320.

\bibitem{levine:book}
M.~Levine:
{\it Mixed motives}, Mathematical Surveys and Monographs, 57,
Amer.\ Math.\ Soc., Providence, RI, 1998.
 
\bibitem{matsumoto}
M.~Matsumoto:
{\it On the Galois image in derivation of $\pi_1$
of the projective line minus three points},
in ``Recent developments in the inverse Galois problem (Seattle, WA, 1993),'' 
Contemp.\ Math.\ 186 (1995), 201--213.

\bibitem{milnor}
J.~Milnor:
{\it Introduction to algebraic $K$-theory}, Annals of Mathematics Studies,
No.~72. Princeton University Press, Princeton, N.J.; University of Tokyo Press,
Tokyo, 1971.

\bibitem{quillen:fg}
D.~Quillen:
{\it Finite generation of the groups $K\sb{i}$ of rings of algebraic integers},
Algebraic $K$-theory, I: Higher $K$-theories (Proc.\ Conf., Battelle Memorial
Inst., Seattle, Wash., 1972), pp. 179--198. Lecture Notes in
Math., Vol. 341, Springer, Berlin, 1973.

\bibitem{quillen:fte}
D.~Quillen:
{\it On the cohomology and $K$-theory of the general linear groups over a
finite field}, Ann.\ of Math.\ 96 (1972), 552--586.

\bibitem{quillen}
D.~Quillen:
{\it Higher algebraic $K$-theory, I}, Algebraic $K$-theory, I: Higher
$K$-theories (Proc.\ Conf., Battelle Memorial Inst., Seattle, Wash., 1972),
pp.~85--147. Lecture Notes in Math., Vol. 341, Springer, Berlin 1973.

\bibitem{Rognes}
J.~Rognes and C.~Weibel:
{\it Two-primary algebraic $K$-theory of rings of integers
in number fields. (Appendix A by Manfred Kolster)},
J.\ Amer.\ Math.\ Soc., 13 (2000), 1--54.

\bibitem{sharifi}
R.~Sharifi:
Letter to Matsumoto, April 28, 2000.

\bibitem{soule}
C.~Soul\'e:
{\it On higher $p$-adic regulators}, Lecture Notes in Math.\ 854 (1981),
372--401.

\bibitem{tate}
J.~Tate:
{\it Relations between $K_2$ and Galois Cohomology},
Invent.\ Math., 30 (1976), 257--274.

\bibitem{tsunogai}
H.~Tsunogai:
{\it On ranks of the stable derivation algebra and
Deligne's problem}, Proc.\ Japan Academy Ser.\ A 73 (1997), 29--31.

\bibitem{voevodsky}
V.~Voevodsky, A.~Suslin, E.~Friedlander:
{\it Cycles, transfers, and motivic homology theories},
Annals of Mathematics Studies, 143, Princeton University Press, 2000.

\bibitem{yoneda2}
N.~Yoneda:
{\it On Ext and exact sequences}, J.\ Fac.\ Sci.\ Univ.\ Tokyo Sect.\ I 8
(1960), 507--576.

\end{thebibliography}
\end{document}